\documentclass[a4paper,10pt]{article}

\usepackage{mymath}
{
	\theoremstyle{plain}
	\newtheorem{assumption}{Assumption}
}

\usepackage[normalem]{ulem}
\ifdefined\C
\renewcommand{\C}{\mathbb{C}}
\else
\newcommand{\C}{\mathbb{C}}
\fi
\title{Efficient function approximation on general bounded domains using splines on a Cartesian grid}
\author{Vincent Copp\'e\footnote{Email: \texttt{vincent.coppe@cs.kuleuven.be}.} \and Daan Huybrechs\footnote{Email: \texttt{daan.huybrechs@cs.kuleuven.be}.}}

\date{KU Leuven \\ Department of Computer Science \\ Celestijnenlaan 200A \\ 3001 Leuven, Belgium \\ \vskip10pt \today}

\newcommand*{\SPAN}{\mathrm{span}}
\newcommand*{\RANK}{\mathrm{rank}}

\newlength{\panelwidth}

\def\bso{p}
\def\osi{q}
\def\bspline{\beta}
\begin{document}
\maketitle
\begin{abstract}
Functions on a bounded domain in scientific computing are often approximated using piecewise polynomial approximations on meshes that adapt to the shape of the geometry. We study the problem of function approximation using splines on a regular but oversampled grid that is defined on a bounding box. This approach allows the use of high order and highly structured splines as a basis for piecewise polynomials. The methodology is analogous to that of Fourier extensions, using Fourier series on a bounding box, which leads to spectral accuracy for smooth functions. However, Fourier extension approximations involve solving a highly ill-conditioned linear system, and this is an expensive step. The computational complexity of recent algorithms is $\mathcal O\left(N\log^2\left(N\right)\right)$ in 1-D and $\mathcal O\left(N^2\log^2\left(N\right)\right)$ in 2-D. We show that, compared to Fourier extension, the compact support of B-splines enables improved complexity for multivariate approximations, namely $\mathcal O(N)$ in 1-D, $\mathcal O\left(N^{3/2}\right)$ in 2-D and more generally $\mathcal O\left(N^{3(d-1)/d}\right)$ in $d$-D with $d>1$. By using a direct sparse QR solver for a related linear system, we also observe that the computational complexity can be nearly linear in practice. This comes at the cost of achieving only algebraic rates of convergence. Our statements are corroborated with numerical experiments and Julia code is available.
\end{abstract}

\noindent \textbf{Keywords} \, function approximation, frames, B-spline, Fourier extension, dual spline basis, oversampling  \\

\noindent \textbf{Mathematics Subject Classification (2010)} \,  65D07, 65D15, 65Y04

\section{Introduction}

It is not straightforward to construct high-order bases for the approximation of functions defined on a domain $\Omega$ with an irregular shape. In contrast, it is conceptually easy to approximate such functions using a basis on a bounding box $\Xi$. Indeed, since a box has Cartesian product structure, one can use the tensor product of any high-order univariate basis and use that to approximate a given function on $\Omega \subset \Xi$.

The choice of classical multivariate Fourier series on $\Xi$ leads to the approximation scheme known as Fourier extension or Fourier continuation~\cite{Boyd2005,Bruno2007,huybrechs2010fourierextension}. The function on $\Omega$ is approximated in a least squares sense and this leads one to solving highly ill-conditioned linear systems. In spite of this apparent conditioning concern, a regularizing singular value decomposition (SVD) -- using a small truncation threshold $\epsilon$ below which singular values are discarded -- can be employed to find accurate and stable approximations~\cite{adcock2014stability,Adcock2019,adcock2020frames}. The direct computation of the SVD has cubic complexity $\mathcal O\left(N^3\right)$, but more recent algorithms for Fourier extensions have been proposed with similar results and much lower computational complexity~\cite{Lyon2011,matthysen2016fast,Matthysen2017}. These methods were generalized to the so-called AZ algorithm for ill-conditioned least squares problems in~\cite{az}.

In this paper, we explore the use of a basis of piecewise polynomials, splines, rather than Fourier series. A similar methodology is used in immersed boundary methods or fictitious domain methods, such as the finite cell method (FCM)~\cite{Parvizian2007,Schillinger2015}. The FCM also defines a simple regular basis of piecewise polynomials on a hypercube that embeds the domain. There are, however, some differences in the FCM approach and the one that we explore. In this paper we are concerned only with the approximation of functions, not the solution of PDEs. It is also worth noting that we make no assumptions on the function in the extended domain. In contrast, fictitious domain methods often extend one or more of the quantities involved, such as the PDE itself, to a formulation on the bounding box, sometimes involving the choice of a parameter or reducing the overall order of the scheme. Our methodology does not depend on the regularity of the function outside $\Omega$, does not rely on an extension of $f$ to the box and readily leads to high-order approximations regardless of the shape of the domain.

Another difference with existing embedding schemes is that we choose to approximate with the unmodified, truncated spline extension frames. This is arguably the simplest and most straightforward approach, which leads to stable and accurate approximation, at the cost of having to solve an ill-conditioned system. Earlier work is mostly focused on designing approximation schemes that circumvent this ill-conditioning. The most prominent example, and most closely related to the current paper, is WEB-splines~\cite{hollig2001b,hollig2003finite}. The method transforms the spline extension frame into a stable spline basis by removing some degrees of freedom. More specifically, the outer B-splines, i.e., the ones that overlap only slightly with the domain (and therefore the sources of ill-conditioning), are adjoined to their inner neighbours. This is achieved in such a way that the resulting basis has full polynomial reproduction on the domain and hence full approximation power. Here, we illustrate that, while seemingly undesirable, the ill-conditioning is harmless if oversampling and regularization is used, see~\cite{Adcock2019,adcock2020frames} for more details.  Furthermore, algorithms are available that are more efficient than a general truncated SVD.

The piecewise polynomials we consider are linear combinations of B-splines. We recall in~\S\ref{s:periodicsplines} the definition of B-splines
and show how bases are formed using translations over a regular, equidistant grid. The dual basis plays a crucial role in our algorithms. 
We recall that the dual basis functions to these B-splines with respect to the inner product of $L^2$ on the real line are not compactly supported, but they do exhibit exponential decay. Proofs of this statement are included in Appendix \ref{app:A} and formalize the ideas in~\cite{Averbuch2016}.
We extend these results to a duality with respect to a discrete inner product in an oversampled equispaced grid, making use of earlier results of~\cite{Unser1991}. 
In the discrete setting as well, the duals are not compactly supported.
We constructively show the existence of compactly supported sequences dual to the B-splines in the oversampled equispaced grid. We extend these results to a periodic setting, which finally forms the basis for the approximation problem formulated in~\S\ref{s:theapproximationproblem}: we restrict (tensor products of) periodic B-spline bases to the domain $\Omega$ at hand, and consider least squares approximants in oversampled grids.

The main part of the paper is the extension in~\S\ref{s:fast} of the above-mentioned AZ algorithm for Fourier extension approximations to these spline-based approximation problems. The AZ algorithm for the solution of a linear system $Ax=b$ relies on knowledge of a $Z$ matrix, such that the matrix $A - AZ^*A$ has numerically low rank. The computational gain of the AZ algorithm results from replacing the full and costly SVD of $A$ by the SVD of (or any other suitable direct solver for) $A - AZ^*A$, and to achieve the latter at a cost that scales with the rank of that matrix rather than with its dimensions. We show that a suitable $Z$ matrix can be derived analytically from the dual bases described earlier for the full periodic basis, with the most efficient results corresponding to the use of a compact discrete dual. 

We introduce two AZ approaches in \S\ref{s:algorithms}. The first is based on a sparse QR factorization of $A - AZ^*A$. This leads to apparent linear complexity in $N$, which however can not be analytically proven. The second consists of a dense QR applied to a reduced linear system whose size is comparable to its rank. The complexity of the latter approach can be fully analysed.
In a third approach, the  sparse QR factorization is directly applied to the system $Ax=b$. This is the most straightforward approach, leads to apparent linear complexity in $N$, and slightly better accuracy. However, the complexity and stability of this approach is not analytically proven.

We illustrate the methods with numerical experiments in~\S\ref{s:experiments}. The results demonstrate a complexity of  $\mathcal O\left(N^{{3(d-1)}/d}\right)$ of the solver in two and three dimensions, as well as algebraic convergence at a rate that depends on the degree of the splines. As an application, we consider the spline-based approximation of elevation data of Belgium and The Netherlands. In this application, the data is given on an equispaced grid, but the grid points are geometrically confined to the shape of both countries. That is precisely the setting of the least squares approximation in this paper.

\section{B-splines and their periodisation}\label{s:periodicsplines}
B-splines can be defined by a recursive relation. The centred B-spline of first order is defined as
\begin{equation*}\label{eq:bsplinedef1}
\bspline^0(t) = \left\{\begin{array}{ll}1,& t\in\left[-\tfrac12,\tfrac12\right)\\0,&\text{otherwise,}\end{array}\right.
\end{equation*}
i.e., it is the characteristic function of the interval $[-\tfrac12,\tfrac12)$. Higher order B-splines are defined by
\begin{equation*}\label{eq:bsplinedef2}
\bspline^\bso(t) = (\bspline^{\bso-1}\star \bspline^0)(t), \quad k>1,
\end{equation*}
where $f\star g$ is the convolution
\begin{equation*}
(f\star g)(t) = \int_{-\infty}^\infty f(\tau)g(t-\tau) d\tau.
\end{equation*}
These B-spline are called centred because their support
\begin{equation}\label{eq:support}
\mysupp \bspline^\bso = \left[-\tfrac{\bso+1}2,\tfrac{\bso+1}2\right]
\end{equation}
is centred around zero. They generate shift-invariant spaces on the real line $S^\bso=\myspan\left\{\Phi^\bso\right\}$, where $\Phi^\bso=\left\{\bspline^\bso(\cdot-k)\right\}_{k\in\Z}$.

Here, splines were defined on the real line. We intend to focus solely on a periodic setting, though other boundaries conditions may also be possible depending on the size of the bounding box. Periodicity enables the use of FFTs in algorithms. Moreover, in the setting of function approximation on a subdomain $\Omega$ of the periodic domain $\Xi$, periodicity on $\Xi$ is not actually a restriction on the functions on $\Omega$ that we can approximate if $\Omega$ is contained in the interior of $\Xi$.

\subsection{Periodisation on $[0,1]$ with continuous and discrete duals}

We define the scaled, normalized and periodised mother function of $\beta^\bso$  as
\begin{equation}\label{eq:periodicspline}
\phi^{\bso}_{N}(t) = \sum_{l \in \Z} N^{1/2} \beta^\bso(N(t -l)).
\end{equation}
The scaling factor $N^{1/2}$ is chosen to ensure that the $L^2$ norms of $\bspline^\bso$ and $\phi^\bso_N$ are equal.
The translates are ${\phi}^{\bso}_{N,k}(t) = {\phi}^\bso_N(t - \frac{k}{N})$
and we consider the sequence of spline bases with $N$ elements
\begin{equation}\label{eq:splinebasis}
\Phi^{\bso}_N = \{\phi^{\bso}_{N,k}(t)\}_{k=0}^{N-1}.
\end{equation}
Note that the infinite sum in~\eqref{eq:periodicspline} reduces to a finite sum in practice owing to the compact support of $\beta^\bso$. By this construction, all functions in the spaces $V^\bso_N :=~\SPAN(\Phi^{\bso}_N)$ are periodic on $[0,1]$.

\paragraph{Continuous duals.} We denote by $\tilde{\Phi}^\bso_N = \{\tilde{\phi}^\bso_{N,k}\}_{k=0}^{N-1}$ the unique dual B-splines in $V^\bso_N$ that are dual to $\Phi^{\bso}_N$ with respect to the standard inner product of $L^2(0,1)$. They can be written as $\tilde{\phi}^{\bso}_{N,k}(t) = \tilde{\phi}^\bso_N(t - \frac{k}{N})$, where the mother function is given in terms of $\Phi^\bso_N$:
\[
\tilde\phi^p_N(t) = \sum_{l\in\Z}N^{1/2}\tilde\beta(N(t-l))= \sum_{k,l\in\Z}c(k)N^{1/2}\beta(N(t-l)-k) = \sum_{k\in\Z}c(k)\phi^p_N\left(t-\tfrac kN\right).
\]
The scaling factor here ensures biorthonormality, i.e, $\langle\phi^\bso_{N,k},\tilde \phi^\bso_{N,l}\rangle_{L^2(0,1)}=\delta_{k,l}$. Because of the periodicity of $\phi^p_N$ the infinite sum can be replaced by a finite one,

\begin{equation}\label{eq:mother_continuous}
\tilde\phi^{\bso}_N(t) = \sum_{k=0}^{N-1} c^{\bso}_N(k)\phi^{\bso}_{N,k}(t),
\end{equation}
with $c^{\bso}_N$ the periodisation of the coefficient sequence in Theorem~\ref{thm:expdecaycontdual}:
\begin{equation}\label{eq:periodiccontcoefficients}
c^{\bso}_N(k) = \sum_{l\in\Z}c^\bso(k+Nl), \quad k=0,\dots, N-1.
\end{equation}

Following theorem shows that these coefficients decay exponentially. 
We use the shortened notation for the inner product: 
\begin{equation}\label{eq:L2innerproduct}
(f,g)_{L^2(\mathbb{R})} =
\int_{-\infty}^\infty f(t) \overline{g(t)}dt.
\end{equation}
\begin{theorem}[Exponential decay of coefficients of the generator dual to the B-spline]\label{thm:expdecaycontdual}
	Let $\tilde\bspline(t)$ be in the form
	\begin{equation}\label{eq:continuousdualexpansion}
	\tilde\bspline(t) = \sum_{k\in\Z}c(k)\bspline^p(t-k),
	\end{equation}
	and let
	\[ \left(\tilde\bspline(\cdot-k),\bspline\right)_{L^2(\mathbb{R})}=\delta_k , \]
	then the sequence $c(k)$ decays exponentially as $|k|\to \infty$.
\end{theorem}

This theorem is an application of a more general statement in Demko's Theorem, see~\cite{demko1977inverses} or \cite[Theorem 4.3]{devore1993constructive}.
Demko's Theorem states that if the inverse of a banded matrix exist, its elements decay exponentially. We show in addition that in the current context the inverse exists. 

\begin{figure}
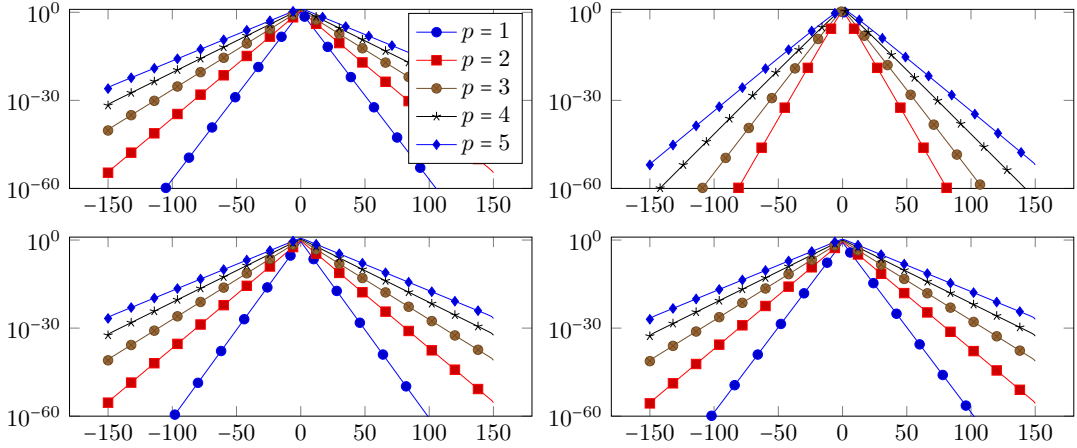

	\setlength{\panelwidth}{.2\columnwidth}
	\centering
	\resizebox{!}{\panelwidth}{\input{P.tikz}}%
	\resizebox{!}{\panelwidth}{\input{Pm1.tikz}}%
	
	\resizebox{!}{\panelwidth}{\input{Pm5.tikz}}%
	\resizebox{!}{\panelwidth}{\input{Pm10.tikz}}%
	\caption{\label{fig:expdecaycoefficients}Size of the dual spline expansion coefficients $c(k)$ (Theorem~\ref{thm:expdecaycontdual}) at top left. Size of the discrete dual spline coefficients $s_\osi(k)$ (Theorem~\ref{thm:expdecaydiscretedual}) for $\osi=1$ (top right), $\osi=5$ (bottom left), and $\osi=10$ (bottom right). The coefficients for $\bso=1,\osi=1$ are not included in the top right panel since they are equal to $\delta_0(k)$. }
\end{figure}

The result is illustrated in the top left of Figure~\ref{fig:expdecaycoefficients} where the exponential decay is clearly visible for several spline orders. For growing order $\bso$ the decay rate is lower. Note that, because B-splines have compact support, exponential decay of the coefficients $c(k)$ in~\eqref{eq:continuousdualexpansion} implies exponential decay of the dual B-spline itself.

Because of the biorthogonality, we can find the coefficients of any $f\in V^\bso_N$ as inner products with the dual basis:
\begin{equation}\label{eq:reconstruction}
	f(t) = \sum_{k \in \mathbb{Z}} ( f , \tilde \phi^\bso_N )_{L^2(\mathbb{R})}\, \phi^\bso_N(t), \qquad f \in V^\bso_N.
\end{equation}

\paragraph{Discrete duals.} We also define an alternative discrete bilinear form
\begin{equation}\label{eq:discreteoversampledinnerproduct}
\langle f, g \rangle_{N,q} = \sum_{k=0}^{Nq-1} f\left(\tfrac k{\osi N}\right)\overline{g\left(\tfrac k{\osi N}\right)}.
\end{equation}
In this discrete setting, we undo the scaling in \eqref{eq:periodicspline} and define $\phi^{\bso,\osi}_{N,k}=1/\sqrt N\phi^{\bso}_{N,k}$, such that the point samples do not depend on $N$. A spline dual with respect to the discrete bilinear form is denoted by $\tilde\phi^{\bso,\osi}_{N,k}(t)=\tilde\phi^{\bso,\osi}_{N}(t-\tfrac kN)$, which means that
\[  
	\left\langle \phi^{\bso,\osi}_{N,k},\tilde \phi^{\bso,\osi}_N\right\rangle_{N,\osi}=\delta_{k}.
\]

Multiple choices for discrete duals to B-splines are possible. The unique dual B-splines in $V^\bso_N$ that are dual to $\Phi^\bso_N$ with respect to the discrete bilinear form have the mother function
\begin{equation}\label{eq:exponentialdiscretedual}
\tilde\phi^{\bso,\osi}_N(t) = \sum_{k=0}^{N-1} c^{\bso,\osi}_N(k)\phi^{\bso,\osi}_{N,k}(t).
\end{equation}
For these discrete dual B-splines, we again see exponential decay if the dual exists, see Theorems~\ref{thm:expdecaydiscretedual} and \ref{thm:completeness} in Appendix~\ref{app:B}. Note that the coefficients $c^{\bso,\osi}_N$ are the periodisation of the coefficient sequences in both theorems,
\begin{equation}\label{eq:periodicdiscretecoefficients}
c^{\bso,\osi}_{N}(k) = \sum_{l\in\Z}s^\bso_\osi(k+Nl), \quad k=0,\dots, N-1.
\end{equation}
Figure~\ref{fig:expdecaycoefficients} illustrates again the effect.

We are interested in duals with compact support, for efficiency reasons. Truncating the exponentially decaying functions using a given threshold $\epsilon>0$ leads to compactly supported duals:
\begin{align}\label{eq:truncateddual}
\tilde\phi^{\bso,\epsilon}_{N}(t) = \left\{\begin{array}{ll}
\tilde\phi^{\bso}_{N}(t)&\text{ if }|t|\leq t^{\bso,\epsilon}_N\\
0&\text{otherwise}
\end{array}\right.,\qquad \tilde\phi^{\bso,\osi,\epsilon}_{ N, k}(t) = \left\{\begin{array}{ll}
\tilde\phi^{\bso,\osi}_{N}(t)&\text{if } |t|\leq t^{\bso,\osi,\epsilon}_N\\
0&\text{otherwise},
\end{array}\right.
\end{align}
with $t^{\bso,\epsilon}_N$ the smallest positive real number $t$ such that {$|\tilde\phi^{\bso}_{N}(t)|<\epsilon$} for $t > t^{\bso}_N$, and $t^{\bso,\osi,\epsilon}_N$ the smallest fraction $t=\tfrac{k}{qN}$ , $k\in\N$ such that {$|\tilde\phi^{\bso,\osi}_{N}(t)|<\epsilon$} for $t > t^{\bso,\osi}_N$.  We define the shifted basis functions $\phi^{\bso,\epsilon}_{N,k}$ and $\phi^{\bso,\osi,\epsilon}_{N,k}$ in terms of these mother functions in the same way as before.  In higher dimensions, the supports of $\tilde\phi^{\mathbf\bso}_{\mathbf N,\mathbf k}$ and $\tilde\phi^{\mathbf\bso,\mathbf\osi,\epsilon}_{\mathbf N,\mathbf k}$ become the smallest hypercubes out of which $|\tilde\phi^{\mathbf\bso}_{\mathbf N,\mathbf k}(\mathbf t)|$ and {$|\tilde\phi^{\mathbf\bso,\mathbf\osi}_{\mathbf N,\mathbf k}(\mathbf t)|$} are smaller then $\epsilon$, respectively.
These duals are biorthogonal up to some numerical error. 

Unfortunately, the decay rate of the coefficients is low. Therefore, the compact supports of $\tilde\phi^{\bso,\epsilon}_{N}$ and $\tilde\phi^{\bso,\osi,\epsilon}_{N}$ are relatively large, which necessitates large $N$ before significant effects are seen in the methods proposed below. 

\begin{figure}
	\centering
	\resizebox{.3\textwidth}{!}{\begin{tikzpicture}
\begin{axis}[legend cell align={left}, mark options={fill opacity={0.2}}]
\addplot+[, samples at={[-10, -9, -8  …  8, 9, 10]}, ycomb, mark={*}]
table[row sep={\\}]
{
	\\
	-10.0  0.0  \\
	-9.0  0.0  \\
	-8.0  0.0  \\
	-7.0  0.0  \\
	-6.0  0.0  \\
	-5.0  0.0  \\
	-4.0  -9.293798490409936e-14  \\
	-3.0  -0.04444444444391689  \\
	-2.0  0.7111111111098443  \\
	-1.0  -3.2777777777761403  \\
	0.0  6.222222222221038  \\
	1.0  -3.2777777777773656  \\
	2.0  0.7111111111110947  \\
	3.0  -0.04444444444455972  \\
	4.0  3.998978879152507e-15  \\
	5.0  0.0  \\
	6.0  0.0  \\
	7.0  0.0  \\
	8.0  0.0  \\
	9.0  0.0  \\
	10.0  0.0  \\
}
;
\addlegendentry {$K=2$}
\addplot+[, samples at={[-10, -9, -8  …  8, 9, 10]}, ycomb, mark={*}]
table[row sep={\\}]
{
	\\
	-10.0  0.0  \\
	-9.0  0.0  \\
	-8.0  0.0  \\
	-7.0  0.0  \\
	-6.0  0.0  \\
	-5.0  0.0  \\
	-4.0  -9.293798490409936e-14  \\
	-3.0  -0.04444444444391689  \\
	-2.0  0.7111111111098443  \\
	-1.0  -3.2777777777761403  \\
	0.0  6.222222222221038  \\
	1.0  -3.2777777777773656  \\
	2.0  0.7111111111110947  \\
	3.0  -0.04444444444455972  \\
	4.0  3.998978879152507e-15  \\
	5.0  0.0  \\
	6.0  0.0  \\
	7.0  0.0  \\
	8.0  0.0  \\
	9.0  0.0  \\
	10.0  0.0  \\
}
;
\addlegendentry {$K=3$}
\addplot+[, samples at={[-10, -9, -8  …  8, 9, 10]}, ycomb, mark={*}]
table[row sep={\\}]
{
	\\
	-10.0  -2.5656260097978976e-15  \\
	-9.0  0.0026629541392250724  \\
	-8.0  -0.04260726622743988  \\
	-7.0  0.1955294873282535  \\
	-6.0  -0.3589994924695739  \\
	-5.0  0.10548450447675753  \\
	-4.0  0.5140108900946416  \\
	-3.0  -0.3417076798736534  \\
	-2.0  -0.8105390918363992  \\
	-1.0  0.4771079025414225  \\
	0.0  1.5181155836534617  \\
	1.0  0.47710790254160657  \\
	2.0  -0.8105390918364976  \\
	3.0  -0.3417076798736484  \\
	4.0  0.5140108900944715  \\
	5.0  0.10548450447703962  \\
	6.0  -0.35899949246967494  \\
	7.0  0.19552948732811806  \\
	8.0  -0.042607266227272636  \\
	9.0  0.002662954139151329  \\
	10.0  1.2431133445204472e-14  \\
}
;
\addlegendentry {$K=10$}
\end{axis}
\end{tikzpicture}}
	\resizebox{.3\textwidth}{!}{\begin{tikzpicture}
\begin{axis}[legend cell align={left}, mark options={fill opacity={0.2}}]
\addplot+[, samples at={[-10, -9, -8  …  8, 9, 10]}, ycomb, mark={*}]
table[row sep={\\}]
{
	\\
	-10.0  0.0  \\
	-9.0  0.0  \\
	-8.0  0.0  \\
	-7.0  0.0  \\
	-6.0  0.0  \\
	-5.0  0.0  \\
	-4.0  -9.293798490409936e-14  \\
	-3.0  -0.04444444444391689  \\
	-2.0  0.7111111111098443  \\
	-1.0  -3.2777777777761403  \\
	0.0  6.222222222221038  \\
	1.0  -3.2777777777773656  \\
	2.0  0.7111111111110947  \\
	3.0  -0.04444444444455972  \\
	4.0  3.998978879152507e-15  \\
	5.0  0.0  \\
	6.0  0.0  \\
	7.0  0.0  \\
	8.0  0.0  \\
	9.0  0.0  \\
	10.0  0.0  \\
}
;
\addlegendentry {$K=2$}
\addplot+[, samples at={[-10, -9, -8  …  8, 9, 10]}, ycomb, mark={*}]
table[row sep={\\}]
{
	\\
	-10.0  0.0  \\
	-9.0  0.0  \\
	-8.0  0.0  \\
	-7.0  -0.0049317482662092245  \\
	-6.0  0.07890797225931266  \\
	-5.0  -0.36147049849005597  \\
	-4.0  0.6545097789859182  \\
	-3.0  -0.15138407556471442  \\
	-2.0  -0.9826361844337895  \\
	-1.0  0.425242746669915  \\
	0.0  1.683524017679219  \\
	1.0  0.42524274666995326  \\
	2.0  -0.9826361844337969  \\
	3.0  -0.15138407556473474  \\
	4.0  0.654509778985907  \\
	5.0  -0.3614704984900317  \\
	6.0  0.07890797225931183  \\
	7.0  -0.004931748266219156  \\
	8.0  0.0  \\
	9.0  0.0  \\
	10.0  0.0  \\
}
;
\addlegendentry {$K=7$}
\addplot+[, samples at={[-10, -9, -8  …  8, 9, 10]}, ycomb, mark={*}]
table[row sep={\\}]
{
	\\
	-10.0  -2.5656260097978976e-15  \\
	-9.0  0.0026629541392250724  \\
	-8.0  -0.04260726622743988  \\
	-7.0  0.1955294873282535  \\
	-6.0  -0.3589994924695739  \\
	-5.0  0.10548450447675753  \\
	-4.0  0.5140108900946416  \\
	-3.0  -0.3417076798736534  \\
	-2.0  -0.8105390918363992  \\
	-1.0  0.4771079025414225  \\
	0.0  1.5181155836534617  \\
	1.0  0.47710790254160657  \\
	2.0  -0.8105390918364976  \\
	3.0  -0.3417076798736484  \\
	4.0  0.5140108900944715  \\
	5.0  0.10548450447703962  \\
	6.0  -0.35899949246967494  \\
	7.0  0.19552948732811806  \\
	8.0  -0.042607266227272636  \\
	9.0  0.002662954139151329  \\
	10.0  1.2431133445204472e-14  \\
}
;
\addlegendentry {$K=10$}
\end{axis}
\end{tikzpicture}}
	\caption{\label{fig:compactdual}The dual signal of $b^1_3$ for $K=2,3,10$ and of $b^4_2$ for $K=2,7,10$.}
\end{figure}
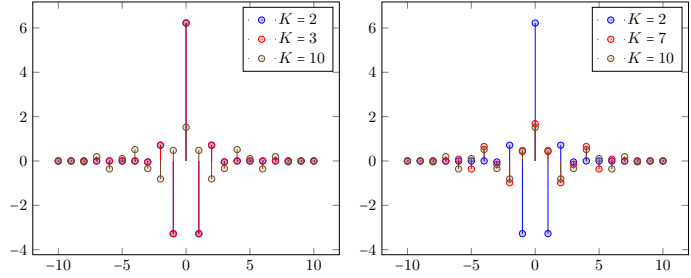
%
%

\paragraph{Compact discrete duals.} Finally, we resort to truly compactly supported discrete duals. They are not in $V^\bso_N$, nor are they unique. 
The discrete duals are known only as sequences $\tilde h^{\bso,\osi}_{N}$, not as continuous functions, but for convenience we associate their values with the evaluations of a dual function in the grid points denoted by
\begin{align}\label{eq:compactsampleddiscrete}
\tilde h^{\bso,\osi}_{N}(k-\osi l) = \tilde\phi^{\bso,\osi}_{N,l}(\tfrac{k}{\osi N}).
\end{align}
As will become clear in the following, their potential values in between grid points does not play a role in our algorithms.

The biorthogonality requirement with respect to the same discrete bilinear form as above,
\[  
\left\langle \tilde\phi^{\bso,\osi}_{N,k},\tilde \phi^{\bso,\osi}_N\right\rangle_{N,\osi}=\delta_{k},
\]
translates into a matrix system  with unknown vector $\tilde h^{\bso,\osi}$, see~\eqref{eq:q-biorthogonality-v2} in Appendix \ref{app:B} for more details. 
Given the dual sequence $\tilde h^{\bso,\osi}$, we can periodise it to $\tilde h^{\bso,\osi}_{N}(k)=\sum_{l\in\Z}\tilde h^{\bso,\osi}(k+N\osi l)$ and obtain the values of $\tilde\phi^{\bso,\osi}_{N}$ in the equispaced grid points. 

For $\osi>1$ there exists a compact sequence $\tilde h^{\bso}_\osi$ dual to $b^\bso_\osi$:
\begin{theorem}\label{thm:compactdual}
	Let $\bso > 0$, $\osi>1$ and let $K$ be a positive integer such that
	\begin{equation}
	K > \frac{p+1}{2}\frac{q}{q-1}-\frac{\osi+1}{q-1}
	\end{equation}
	for $\bso$ odd, and
	\begin{equation}
	K > \left\lceil\frac{p+1}{2}\frac{q}{q-1}-\frac{\osi+1}{q-1}\right\rceil
	\end{equation}
	for $\bso$ even. Then, there exists at least
	one compact sequence $\tilde h^{\bso}_\osi$ dual to $b^\bso_\osi$
	 with support $[-K,K]$.
\end{theorem}

The proof of the theorem is presented in Appendix~\ref{app:B},~\S\ref{ss:compactdual}. Figure~\ref{fig:compactdual} shows some of these duals with different support for $b^1_3$ and $b^3_2$. The results are reminiscent of quasi-interpolation and corresponding families of compactly supported biorthogonal functionals~\cite{guide1978splines,schumaker1981spline}, one difference being our use of an oversampled grid.

In the remainder of the text, we will drop the superscript $\bso$ to make notation lighter.

\begin{remark}\label{rem:dualsize}
There are many discrete compact duals, and Theorem~\ref{thm:compactdual} provides only the minimal support such that existence is guaranteed. As can be seen in Figure~\ref{fig:compactdual}, the size of the dual sequence depends on the support length, with sequences of larger support having smaller elements. This has a beneficial effect on the accuracy of the reconstruction. We will briefly illustrate this effect in \S\ref{s:experiments}.
\end{remark}

\subsection{Computing the dual spline bases}\label{ss:efficient_fft}

The coefficients $c_{N}$ of the continuous dual spline mother function~\eqref{eq:mother_continuous} may be found by inverting the Gram matrix of $\Phi_{N}$, which is given by
\begin{equation}\label{eq:primal_gram}
G_{N}(k,l)=(\phi_{N,k},\phi_{N,l})_{L^2(0,1)}.
\end{equation}
This standard observation follows from the biorthogonality property:
\begin{equation*}
\delta_{k,l}=(\phi_{N,k},\tilde\phi_{N,l})_{L^2(0,1)}  = \sum_{i=1,\dots,N} G_{N}(k,i)c_{N}(i-l),
\end{equation*}
which in matrix notation corresponds to $I=G_{N}M$ with $M(k,l)=c_{N}(k-l)$. The inverse matrix $M=(G_{N})^{-1}$ is in fact the Gram matrix of $\tilde\Phi_{N}$ with elements $\tilde G_{N}(k,l)=(\tilde\phi_{N,k},\tilde\phi_{N,l})_{L^2(0,1)}$, since
\begin{equation*}
(\tilde\phi_{N,k},\tilde\phi_{N,l})_{L^2(0,1)}  = \sum_{i,j=1,\dots,N} c_{N}(i-k)G_{N}(i,j)c_{N}(j-l) = (MM^{-1}M)(k,l) = M(k,l).
\end{equation*}

The Gram matrix and its inverse are circulant because of the periodicity and shift-invariance of $\Phi_{N}$. Hence, they are diagonalized by the FFT and have a fast matrix-vector multiply based on FFTs~\cite{Davis1979,Chen1987a}. Similar reasonings exist for the coefficients in~\eqref{eq:periodicdiscretecoefficients} based on the discrete bilinear form~\eqref{eq:discreteoversampledinnerproduct}. The coefficients in~\eqref{eq:periodiccontcoefficients} and~\eqref{eq:periodicdiscretecoefficients} can thus be computed in $\mathcal O(N\log(N))$ operations. 

The number of operations in the computation of $\tilde h^{\bso,\osi}_{N}$ in~\eqref{eq:compactsampleddiscrete} is independent of $N$ 
as it requires solving a system with a size independent on $N$, see system~\eqref{eq:compact_system} in Appendix~\ref{app:B}.
The size of this system depends only on $\bso$ and $\osi$. 


\section{The approximation scheme}\label{s:theapproximationproblem}

We are interested in the approximation of a function $f$ on a compact domain $\Omega$, where $\Omega$ may have an arbitrary shape. Without loss of generality, we can scale $\Omega$ such that $\Omega\subset[0,1]^d$. Here, $d$ is the dimension of $\Omega$. Whereas Fourier extensions are based on tensor product Fourier series, in this section we study the use of (tensor products of) periodic and shifted B-splines as described in~\S\ref{s:periodicsplines}.

\subsection{Discretization of the approximation problem}

In one dimension, the approximation to $f \in L^2(\Omega)$ with $\Omega \subset [0,1]$ we are interested in is the best approximation
\begin{equation}\label{eq:approximation}
y = \argmin_{a\in\C^N}\left\|f-\sum_{k=0}^{N-1}a(k)\phi_{N,k}\right\|_{L^2(\Omega)},
\end{equation}
{where $\phi_{N,k}$ are the spline basis functions defined by~\eqref{eq:splinebasis}.}

The discretization of this problem leads to the system of equations $A_{N}x=b_{N}$, with matrix and vector entries given in terms of inner products on $\Omega$:
\begin{align}
	A_{N}(k,l) &= (\phi_{N,k},\phi_{N,l})_{L^2(\Omega)}
	\label{eq:galerkinmatrix},\\
	b_{N}(k) &=  (f, \phi_{N,k})_{L^2(\Omega)},\nonumber
\end{align}
where $k,l=0,\ldots,N-1$.

Alternatively, we consider a discrete approximation of~\eqref{eq:approximation}. To that end, we first choose a set of points
\begin{equation}\label{eq:points_univariate}
P^\osi_N := \left\{ \tfrac{k}{q N} \, \big| \, \tfrac{k}{q N} \in \Omega, k\in\Z \right\},
\end{equation}
which arises from the restriction of an equispaced grid with $L=qN$ points on $[0,1]$ to $\Omega$. For convenience we sort the points by size and denote them by $P^\osi_N = \{ t_m \}_{m=0}^{M-1}$ where $M = \# P^\osi_N$. Using these points we consider a discrete least squares approximation
\begin{equation*}\label{eq:approximation_discrete}
y = \argmin_{a\in\C^N} \sum_{m=1}^M \left| f(t_m) - \sum_{k=0}^{N-1}a(k)\phi_{N,k}(t_m)\right|^2.
\end{equation*}
This leads to having to solve the linear system $A^{\osi}_{N} x=b^{\osi}_{N}$ with
\begin{align}
	A^{\osi}_{N}(m,l) &=
	 \phi^\osi_{N,l}\left(t_m\right)\label{eq:collocationmatrix}\\
	b^{\osi}_{N}(m) &= 
	f\left(t_m\right),\nonumber
\end{align}
where $l=0,\ldots,N-1$ and $m=0,\ldots,M-1$. 

In the remainder of the text we assume a linear oversampling, i.e., $M = \mathcal O(N)$. This implies that we determine a number $\osi$ and corresponding factor $L = qN$, such that the set \eqref{eq:points_univariate} has a number of points that scales linearly with $N$, with a proportionality constant greater than $1$. This linear scaling is assumed in the computational complexity statements later on.

\subsection{Multivariate approximation}

We extend the notation above to two and more dimensions. We define a tensor product basis for $[0,1]^d$ using $d$ univariate B-spline bases $\Phi^{\bso_i}_{N_i}$ by
\begin{equation*}
\Phi^{\mathbf\bso}_{\mathbf N} = \Phi^{\bso_1}_{N_1}\otimes\Phi^{\bso_2}_{N_2}\otimes\dots\otimes\Phi^{\bso_d}_{N_d}=\{\phi^{\mathbf\bso}_{\mathbf N, \mathbf k}(\mathbf t)\}_{\mathbf k\in I_{\mathbf N}}
\end{equation*}
with the multi-index set
\begin{equation*}
	I_{\mathbf N} = \{(n_1,\dots,n_d) \, | \, n_i=0,\dots,N_d-1, i=1\dots,d\}.
\end{equation*}
Similarly, the dual bases are tensor products of their univariate versions,
\begin{align*}
	{\tilde\Phi}^{\mathbf\bso}_{\mathbf N}&= \tilde{\Phi}^{\bso_1}_{N_1}\otimes\tilde{\Phi}^{\bso_2}_{N_2}\otimes\dots\otimes\tilde{\Phi}^{\bso_d}_{N_d}=\{\tilde\phi^{\mathbf\bso}_{\mathbf N,\mathbf k}(\mathbf t)\}_{\mathbf k\in I_{\mathbf N}}\\
	{\tilde\Phi}^{\mathbf\bso,\mathbf\osi}_{\mathbf N}&= \tilde{\Phi}^{\bso_1,\osi_1}_{N_1}\otimes\tilde{\Phi}^{\bso_2,\osi_2}_{N_2}\otimes\dots\otimes\tilde{\Phi}^{\bso_d,\osi_d}_{N_d}=\{\tilde\phi^{\mathbf\bso,\mathbf\osi}_{\mathbf N,\mathbf k}(\mathbf t)\}_{\mathbf k\in I_{\mathbf N}}.
\end{align*}

The continuous multivariate approximation problem for $\Omega \subset [0,1]^d$ leads to the system $A_{\mathbf N}x=b_{\mathbf N}$ where $A_{\mathbf N}\in\R^{N\times N}$ (with $N = \Pi_{i=1,\dots,d}N_i$) and
\begin{align*}
	A_{\mathbf N}(\mathbf k,\mathbf l) &= (\phi_{\mathbf N,\mathbf k},\phi_{\mathbf N,\mathbf l})_{L^2(\Omega)}, \\
	b_{\mathbf N}(\mathbf k) &= (f, \phi_{\mathbf N,\mathbf k})_{L^2(\Omega)},
\end{align*}
with $\mathbf k,\mathbf l\in I_{\mathbf N}$.

For the discrete approach, we first determine $M$ points $P^{\mathbf\osi,\Omega}_{\mathbf N} = \{\mathbf t_m\}_{m=0}^{M-1}$ that are in the intersection of $\Omega$ and the Cartesian grid
\begin{equation*}
P^{\mathbf\osi}_{\mathbf N} = \left\{\left.\left(\tfrac{k_1}{\osi_1 N_1},\dots,\tfrac{k_d}{\osi_d N_d}\right)\, \right| \, k_i=0,\dots,N_i\osi_i-1,\quad i=1,\dots, d \right\}.
\end{equation*}
Next, we construct the linear system $A^{\mathbf\osi}_{\mathbf N} x = b^{\mathbf\osi}_{\mathbf N}$ with $A^{\mathbf\osi}_{\mathbf N}\in\R^{M\times N}$ and
\begin{align*}
A^{\mathbf\osi}_{\mathbf N}(m,\mathbf l) &= 
\phi^{\mathbf\osi}_{\mathbf N,\mathbf l}(\mathbf t_m)\\
b^{\mathbf\osi}_{\mathbf N}(m) &= 
f(\mathbf t_m)
\end{align*}
where $m=I_{M}$ and $\mathbf l\in I_{\mathbf N}$. 

The continuous approximation problem leads to the practical difficulty of having to evaluate integrals on $\Omega$, in order to evaluate the inner products. In contrast, for the discrete problem it suffices to be able to determine whether or not a point on the Cartesian grid above belongs to $\Omega$. In addition, it must be possible to sample $f$ at such points.

Here, too, we will use a linear sampling rate $M = \mathcal O(N)$, with proportionality constant greater than $1$ to ensure oversampling.

\subsection{Numerical solution methods}

Though the basis functions are linearly independent on $[0,1]^d$, this need not be the case on $\Omega$. As a result, we can expect that the matrices $A_{\mathbf N}$ and $A^{\mathbf\osi}_{\mathbf N}$ are ill-conditioned, depending on $\Omega$ and on the choice of the parameters ${\mathbf N}$, ${\mathbf p}$ and ${\mathbf q}$.

Following the recommendation of~\cite{Adcock2019} for the approximation of functions using redundant sets, we solve the linear systems using a SVD of $A$, with all singular values below a threshold $\epsilon$ discarded. This method is straightforward to implement but it has a computational complexity of ${\mathcal O}(N^3)$ and we set out to reduce this complexity.

\section{Fast approximation using the AZ algorithm}\label{s:fast}

\subsection{The AZ algorithm for matrices with a plunge region}

\begin{algorithm}[t]
	\caption{The AZ algorithm~\cite{az}}\label{alg:AZ}
	{\bf Input:} $A,Z \in \C^{M\times N}$, $b\in\C^M$ \\
	{\bf Output:} $x\in\C^N$ such that $Ax \approx b$
	\begin{algorithmic}[1]
		\State Solve $(I-AZ^*)Ax_1 = (I-AZ^*)b$
		\State $x_2 \gets Z^*(b-Ax_1)$
		\State $x \gets x_1 + x_2$
	\end{algorithmic}
\end{algorithm}

The AZ algorithm was introduced in~\cite{az} as a generic way to solve least squares systems with system matrix $A$ for which an approximate inverse $Z^*$ is known. The AZ algorithm is shown in Algorithm~\ref{alg:AZ}. The notion of `approximate inverse' in this context means the following. The matrix $Z$ can in principle be chosen freely, but the goal is to ensure that $A - AZ^*A$ has low (numerical) rank. Note that $A - AZ^*A$ appears in step 1 of the algorithm with a direct solve, which can be implemented efficiently using a randomized linear algebra solver if the matrix indeed has low rank. The goal is achieved if $Z^*$ is an inverse for a large subspace $Y$ of the column space of $A$, i.e., if for $A \in \mathbb{C}^{M \times N}$ one has that $Z^*A y \approx y$ for $y \in Y \subset \mathbb{C}^N$. It is observed in~\cite{az} with a variety of examples that approximation problems with known efficient solvers can be leveraged to more general problems for which the original solver is no longer exact, but does remain well suited as a $Z$ matrix. The AZ algorithm is most efficient if both $A$ and $Z$ have a fast matrix-vector multiply.

A large class of matrices for which the AZ algorithm applies involve matrices with a \emph{plunge region}. The plunge region refers to a sudden drop of the singular values from a region in which they are ${\mathcal O}(1)$ down to zero. The plunge region has received the most study in the context of Fourier extensions and prolate spheroidal wave functions~\cite{slepian1978,edelman1999}. In these applications, $Z$ is chosen to \emph{isolate the plunge region}: the rank of $A - AZ^*A$ scales with the width of the plunge region, rather than with the dimension of $A$. If the non-decaying singular values of $A$ (before the onset of the plunge region) are approximately $1$, as is the case when using normalized Fourier series in a Fourier extension problem, then one can simply choose $Z=A$~\cite{matthysen2016fast,Matthysen2017}. Matrix-vector products involving either $A$ or $Z$ can be expedited using the FFT.

In the remainder of this section, we set out to show that the system matrices of the spline extension approximation scheme exhibit a plunge region, and that a suitable $Z$ matrix can be found analytically or numerically from the dual basis on the bounding box. Fast matrix-vector multiplies are obtained either by using sparsity or by exploiting circulant structure using the FFT.

\subsection{The choice of $Z$ based on duality}

We first detail the choice of $Z$ for the continuous and discrete spline extension approximation schemes. In both cases, it is based on the dual of the original basis on the full domain $[0,1]^d$. The proof that these choices lead to low rank of $A Z^* A - A$ follows in the next section.

In the continuous approach, the system matrix $A$ is $A_N$~\eqref{eq:galerkinmatrix}.
We simply choose $Z_N = \tilde G_N$, where $\tilde G_N$ is the Gram matrix of the dual basis $\tilde{\Phi}_N$ as defined in~\S\ref{ss:efficient_fft}. Thus, the elements of $Z$ are inner products on the full domain $[0,1]$ rather than on the subset $\Omega$:
\begin{equation}\label{eq:galerkingZ}
	Z_{N}(k,l) = \tilde G_{N}(k,l)= (\tilde\phi_{N,k},\tilde\phi_{N,l})_{L^2(0,1)}, \qquad k,l=0,\ldots,N-1.
\end{equation}
We have by construction of the dual that $Z^*_N=(\tilde G^{\bso}_{N})^*$ is the inverse of the Gram matrix $G_N$ of the primal basis $\Phi_N$ on $[0,1]$, recall~\eqref{eq:primal_gram}. As we will show, owing to the compact support of B-splines, this property also ensures that $Z^*$ is an approximate inverse for $A_N$ suitable for the AZ algorithm. Indeed, this is the most crucial observation of the current paper.

In higher dimensions the $Z$ matrix is denoted $Z_{\mathbf N}$. It is the Kronecker product of $Z_{N_i}$ for $i=1,\dots,N$, or
\begin{equation*}
	Z_{\mathbf N}(\mathbf k,\mathbf l) = (\tilde\phi_{\mathbf N,\mathbf k},\tilde\phi_{\mathbf N,\mathbf l})_{L^2([0,1]^d)},
\end{equation*}
where $\mathbf k,\mathbf l\in I_{\mathbf N}$. Note that $A_{\mathbf N}$ itself does not have such Kronecker structure, unless $\Omega$ happens to be a Cartesian product domain.

In the discrete approach, the system matrix $A$ is $A^{\osi}_{N}$, given by~\eqref{eq:collocationmatrix}.
Here, a fast matrix-vector product is available since  $A^{\osi}_{N}$ contains only $\mathcal O(N)$ non-zero elements. The $Z$ matrix that leads to a low-rank $A-AZ^*A$ is based on sampling the discrete dual functions,
\begin{equation}\label{eq:collocationZ}
	Z^{\osi}_{N}(m, l)= 
	\tilde\phi^{\osi}_{N,l}\left(t_m\right),
\end{equation}
where $t_m$ are the points of the set $P^\osi_N$, recall~\eqref{eq:points_univariate}.
Using the first discrete dual~\eqref{eq:exponentialdiscretedual}, $Z$ is a full matrix, but it can be applied fast since it is the submatrix of a circulant matrix. Using the truncated discrete dual~\eqref{eq:truncateddual}, $Z$ is banded, but it typically has unpractically large bands. Finally, using the compact discrete dual~\eqref{eq:compactsampleddiscrete}, $Z$ is a banded matrix with a relatively small bandwidth.

In higher dimensions we choose $Z$ as $Z^{\mathbf\osi}_{\mathbf N}$ with the analogous definition
\begin{equation*}
	Z^{\mathbf\osi}_{\mathbf N}(m,\mathbf l) =
	\tilde\phi^{\mathbf\osi}_{\mathbf N,\mathbf l}(\mathbf t_m),
\end{equation*}
where $m\in I_M$, $\mathbf t_m\in P^{\mathbf\osi}_{\mathbf N}$ and $\mathbf l\in I_{\mathbf N}$.

If the notation is clear from the context, we will use $A$ and $Z$ to denote either the combination $A_{\mathbf N}$ and $Z_{\mathbf N}$, or the combination $A^{\mathbf\osi}_{\mathbf N}$ and $Z^{\mathbf\osi}_{\mathbf N}$.

\section{Two AZ algorithms for spline extension}\label{s:algorithms}

Earlier algorithms based on the AZ algorithm focused on the numerical low rank of $A-AZ^*A$. 
However, in the case of spline extension, as we will show in Theorem~\ref{thm:lowrank1}, the matrix $A-AZ^*A$ is highly sparse and even contains a lot of zero columns and zero rows. The algorithms we introduce here take advantage of these properties.

The first algorithm, the \emph{sparse AZ algorithm}, follows the same structure as Algorithm~\ref{alg:AZ} but builds a sparse representation of $A-AZ^*A$. It solves the first step using a direct least squares solver that relies on a rank-revealing sparse QR factorization~\cite{Davis2009} of $A-AZ^*A$. This routine is available, e.g.,  in Matlab as \texttt{spqr\_solve}.

It is difficult to make definitive statements regarding the time complexity and stability of this solver based on a sparse QR factorization. Hence, we also consider an algorithm that uses a traditional pivoted QR factorization, which we denote the \emph{reduced AZ algorithm} (Algorithm~\ref{alg:reducedAZ}). Knowing the indices of the non-zero rows and columns, the factorization in step 1 of the AZ algorithm can be performed on a submatrix  $R(A-AZ^*A)E$ of reduced size. This submatrix contains only the non-zero rows ($R$) and non-zero columns ($E$).

Of course the matrix $A$ itself is also highly sparse. In the experimental results, we will compare against a \emph{direct sparse QR solver} applied to the system $Ax=b$ itself. This is a simpler and more straightforward approach which, as we shall see, is surprisingly competitive with the more advanced algorithms. However, as with the sparse AZ algorithm, its computational complexity can not be proven.


\subsection{The structure of $A-AZ^*A$}\label{ss:lowrank}

The efficiency gain of the AZ algorithm arises from the low rank of $A-AZ^*A$. Furthermore, the proposed algorithms depend on the sparsity structure of $A-AZ^*A$. Hence we start by analysing this matrix. The main result is formulated in Theorem~\ref{thm:lowrank1}. Note that the statement of the theorem is independent of the duals proposed in~\S\ref{s:periodicsplines}. 

First, we define the index set of spline basis functions whose support overlaps with the boundary of a domain $S$:
\begin{equation}
 \label{eq:boundaryset}
 K^{\mathbf\bso}_{\mathbf N}(S) = \{ {\mathbf l} \in I_{\mathbf N} \, | \, \, \mysupp \phi^{\mathbf\bso}_{\mathbf N,\mathbf l} \cap S \neq \varnothing\text{ and }\mysupp \phi^{\mathbf\bso}_{\mathbf N,\mathbf l} \cap {S^c} \neq \varnothing \}.
\end{equation}
with $S^c = [0,1]^d \setminus S$. We use the notation $K_{\mathbf N}(S)$ to emphasize here that the quantity depends on the domain $S$, but for notational convenience we may omit the domain further on.

In order to analyse the rank of $A-AZ^*A$ in the discrete setting, we also introduce a notion of \emph{discrete support}. The discrete support is the intersection of the support of a basis function with the grid points in $P^{\mathbf\osi}_{\mathbf N}$:
\begin{equation*}
\mysupp_{\mathbf\osi} \phi^{\mathbf\osi}_{\mathbf N,\mathbf l} = \mysupp\phi^{\mathbf\osi}_{\mathbf N,\mathbf l} \cap P^{\mathbf\osi}_{\mathbf N}.
\end{equation*}

Generalizing~\eqref{eq:boundaryset}, we define the set of spline basis functions whose discrete support overlaps with the boundary of a domain $S$ as:
\begin{equation}\label{eq:boundaryset_discrete}
K_{\mathbf N}^{\mathbf\bso,\mathbf\osi}(S) = \{ \mathbf l \in I_{\mathbf N}, \, | \, \mysupp_{\mathbf\osi} \phi^{\mathbf\bso,\mathbf\osi}_{\mathbf N,\mathbf l} \cap S \neq \varnothing\text{ and }\mysupp_{\mathbf\osi} \phi^{\mathbf\bso,\mathbf\osi}_{\mathbf N,\mathbf l} \cap {S^c} \neq \varnothing \}.
\end{equation}

Now we are ready to state the main theorem on the structure of $A-AZ^*A$:

\begin{theorem}\label{thm:lowrank1}
        For the AZ pair $(A^{\mathbf\bso}_{\mathbf N},Z^{\mathbf\bso}_{\mathbf N})$, and some $\delta>0$,
        the matrix $A-AZ^*A$ has 
        \begin{enumerate}
        	\item at most $\# K^{\mathbf\bso}_{\mathbf N}(\Omega)$ non-zero columns
        	\item $\RANK(A-AZ^*A) \leq \# K^{\mathbf\bso}_{\mathbf N}(\Omega)$
        	\item $\mathcal O(\#  K^{\mathbf\bso}_{\mathbf N}(\Omega))$ rows with elements larger than 
        	$\delta$ 
        	in absolute value and
        	\item $\mathcal O(\#  K^{\mathbf\bso}_{\mathbf N}(\Omega))$ elements larger than 
        	$\delta$ 
        	in absolute value.
        \end{enumerate} 
	For the AZ pair  $(A^{\mathbf\osi}_{\mathbf N},Z^{\mathbf\osi}_{\mathbf N})$, and some $\delta>0$,
	the matrix $A-AZ^*A$ has  
	\begin{enumerate}
		\item at most $\# K_{\mathbf N}^{\mathbf\bso,\mathbf\osi}(\Omega)$ non-zero columns
		\item $\RANK(A-AZ^*A) \leq \# K_{\mathbf N}^{\mathbf\bso,\mathbf\osi}(\Omega)$
		\item $\mathcal O(\# K_{\mathbf N}^{\mathbf\bso,\mathbf\osi}(\Omega))$ rows with elements larger than 
		$\delta$ 
		in absolute value and
		\item $\mathcal O(\# K_{\mathbf N}^{\mathbf\bso,\mathbf\osi}(\Omega))$ elements larger than 
		$\delta$ 
		in absolute value.
	\end{enumerate} 
The constant in the big $\mathcal O$
is independent of $N$.
\end{theorem}

\begin{proof}[Proof of Theorem~\ref{thm:lowrank1}]
	The proofs of the two parts in Theorem~\ref{thm:lowrank1} follow the same reasoning. 
	We consider the continuous result first and then indicate where the proof differs for the discrete case.
	
	\emph{Proof of statements 1. and 2.}
	For the continuous approximation problem, the AZ pair is  $(A_{\mathbf N},Z_{\mathbf N})$.
	We have by construction that
	\begin{equation*}
	f= \sum_{{\mathbf i} \in I_{\mathbf N}} (f,\tilde{\phi}_{\mathbf{N},\mathbf{i}})_{L^2([-1,1]^d)} \phi_{\mathbf{N},\mathbf{i}}
	\end{equation*}
	for all $f \in~\SPAN \, \Phi_{\mathbf{N}}$. Choosing $f={\phi}_{\mathbf{N},\mathbf{k}}$ and applying $(f,g)_{L^2(\Omega)}$ with $\tilde {\phi}_{\mathbf{N},\mathbf{l}}$ shows that
	\begin{align}\label{eq:contZA}
	\left(Z^*A\right)(\mathbf k,\mathbf l) = \sum_{\mathbf i\in I_{\mathbf N}}(\tilde\phi_{\mathbf N,\mathbf k},\tilde\phi_{\mathbf N,\mathbf i})_{L^2([0,1]^d)} (\phi_{\mathbf N,\mathbf l},\phi_{\mathbf N,\mathbf i})_{L^2(\Omega)} = (\phi_{\mathbf N,\mathbf l},\tilde\phi_{\mathbf N,\mathbf k})_{L^2(\Omega)}.
	\end{align}
	Note that $(\phi_{\mathbf N,\mathbf l},\tilde\phi_{\mathbf N,\mathbf k})_{L^2(\Omega)} = (\phi_{\mathbf N,\mathbf l},\tilde\phi_{\mathbf N,\mathbf k})_{L^2([0,1]^d)} = \delta_{\mathbf k, \mathbf l}$ whenever $\mysupp \phi_{\mathbf N,\mathbf l} \subset \Omega$. Similarly, the inner product evaluates to zero if $\phi_{\mathbf N,\mathbf l}$ is supported outside of $\Omega$.
	
	Let ${\mathbb I}_{\mathbf N}$ be the Kronecker product of the identity matrices of size $N_i$ for $i=1,\dots,d$. {Then, the above reasoning demonstrates sparsity of the matrix}
	\begin{align*}
	\left({\mathbb I}_{\mathbf N} - Z^*A\right)(\mathbf k,\mathbf l) = \left\{
	\begin{array}{ll}
	\delta_{\mathbf k,\mathbf l} & \text{if }\mysupp\phi_{\mathbf N,\mathbf l}\subset[0,1]^d\backslash\Omega,\\
	0 & \text{if }\mysupp\phi_{\mathbf N,\mathbf l}\subset\Omega,\\
	\delta_{\mathbf k,\mathbf l}-(\phi_{\mathbf N,\mathbf l},\tilde\phi_{\mathbf N,\mathbf k})_{L^2(\Omega)} &\text{otherwise}.
	\end{array}\right.
	\end{align*}
	
	It follows that $I-Z^*A$ has a zero column with index $\mathbf{l}$ if the basis function $\phi_{\mathbf N,\mathbf l}$ is contained within $\Omega$, or a unit-vector column if $\phi_{\mathbf N,\mathbf l}$ is contained in the complement of $\Omega$. Another multiplication by $A$ leads to even more sparsity:
	\begin{align}\label{eq:plungestructure2}
	\left( A (I - Z^*A) \right)(\mathbf k,\mathbf l) = \left\{
	\begin{array}{ll}
	0 & \text{if }\mysupp\phi_{\mathbf N,\mathbf l}\subset[0,1]^d\backslash\Omega,\\
	0 & \text{if }\mysupp\phi_{\mathbf N,\mathbf l}\subset\Omega,\\
	a_{\mathbf k,\mathbf l} &\text{otherwise},
	\end{array}\right.
	\end{align}
	where the values $a_{\mathbf k,\mathbf l}$ may or may not be zero. The first line follows because $A(\mathbf k,\mathbf l)=0$ if $\mysupp\phi_{\mathbf N,\mathbf l}\subset[0,1]^d\backslash\Omega$ for $\mathbf k\in I_{\mathbf N}$, hence the multiplication by the unit vector in column $\mathbf l$ of $(I-Z^*A)$ evaluates to zero.
	
	The number of columns in $A(I-Z^*A)$ that are not identically zero is bounded by the number of basis functions in the set $K_{\mathbf N}(\Omega)$, because only these basis functions do not match either of the first two conditions in~\eqref{eq:plungestructure2}.
	
	This concludes the proof on the statement of the number of non-zero columns. The bound on the rank follows immediately.

	For the discrete case, it suffices to replace the AZ pair by $(A^{\mathbf\osi}_{\mathbf N},Z^{\mathbf\osi}_{\mathbf N})$, to use the discrete bilinear forms and supports: $\langle f,g\rangle_{N,\mathbf\osi,\Omega} $, $\langle f,g\rangle_{N,\mathbf\osi} $ and $\mysupp_{\mathbf \osi}$, instead of $(f,g)_{L^2(\Omega)}$, $(f,g)_{[-1,1]^d}$; and to replace $K_{\mathbf N}(\Omega)$ by $K_{\mathbf N}^{\mathbf\osi}(\Omega)$. The same substitutions are largely sufficient for the second half of the proof, with additional remarks explicitly noted.

	\emph{Proof of statements 3. and 4.}
	Here we introduce 
	\begin{align*}
	Z^{\mathbf\bso,\epsilon}_{\mathbf N}(\mathbf k,\mathbf l) = (\tilde\phi^{\mathbf\bso,\epsilon}_{\mathbf N,\mathbf k},\tilde\phi^{\mathbf\bso,\epsilon}_{\mathbf N,\mathbf l})_{L^2([0,1]^d)}, 
	\end{align*}
	which uses the truncated versions in \eqref{eq:truncateddual} of the dual function used to construct $Z^{\mathbf p}_{\mathbf N}$. For the discrete case,  this would be $Z^{\mathbf\bso,\mathbf\osi,\epsilon}_{\mathbf N}(m,\mathbf l) =
	 \tilde\phi^{\mathbf\bso,\mathbf\osi,\epsilon}_{\mathbf N,\mathbf l}(\mathbf t_m)$.
	
	The error made in $Z$ by replacing infinite supports with finite ones can be bounded elementwise by:
	\[
	|(Z_{\mathbf N}-Z^{\epsilon}_{\mathbf N})(\mathbf k,\mathbf l)| = \left| (\tilde\phi_{\mathbf N,\mathbf k},\tilde\phi_{\mathbf N,\mathbf l})_{L^2\left([0,1]^d \setminus \left(\mysupp \tilde\phi^{\epsilon}_{\mathbf N,\mathbf k}\cup\,  \mysupp \tilde\phi^{\epsilon}_{\mathbf N,\mathbf l}\right)\right)} \right| <\epsilon.
	\]
	This is because $\tilde \phi_{\mathbf N,\mathbf k}$ is smaller than $\epsilon$ by construction away from the support of $\tilde\phi^{\epsilon}_{\mathbf N,\mathbf k}$, and the inner product above amounts to the integral of a quantity less than $\epsilon$ over a subset of the domain $[0,1]^d$ which has volume $1$. The elementwise error bound of $\epsilon$ follows in the discrete case from the definition~\eqref{eq:truncateddual}.
	
	Using this elementwise bound, the entries of $A-AZ^*A$ can in turn be bounded by
	\[
	\left| \left(A (I - Z^*A) \right)-\left( A (I - (Z^{\epsilon})^*A) \right)(\mathbf k,\mathbf l) \right| < \epsilon\left|(A  1_{\mathbf N} A)(\mathbf k,\mathbf l)\right|,
	\]
	where $1_{\mathbf N}$ is a matrix with all entries equal to $1$. For the continuous case the entries of $A=A^{\mathbf\bso}_{\mathbf N}$ are all positive and that $A$ is symmetric, we have
	\[
	\left|(A^{\mathbf\bso}_{\mathbf N}  1_{\mathbf N} A^{\mathbf\bso}_{\mathbf N})(\mathbf k,\mathbf l)\right| =  \sum_{\mathbf i \in I_{\mathbf N}}  A^{\mathbf\bso}_{\mathbf N}(\mathbf k,\mathbf i) \sum_{\mathbf j \in I_{\mathbf N}} A^{\mathbf\bso}_{\mathbf N}(\mathbf l,\mathbf j)\leq M^2\Pi_{i=1,\dots,d}{(2\bso_i+1)^2}.
	\]
	Here, $M$ is the largest element of $A$, and $2p_i+1$ is the number of overlapping neighbouring spline functions with $\phi_{\mathbf N,\mathbf k}$ in dimension $i$. The product of these factors in all dimensions is the maximal number of non-zero entries in row $\mathbf k$ of $A$. Since the splines are strictly decreasing away from their centre, the maximal entry of $A$ is $M = \prod_{i=1}^d \|{\bspline^{\bso_i}}\|^2$, which is achieved by any diagonal entry of $A$ associated with a spline that is completely supported within $\Omega$. We denote the upper bound in the expression above as $C$. It is independent on $N$.
	
	The elements of $A^{\mathbf\osi}_{\mathbf N}1_{N\times M}A^{\mathbf\osi}_{\mathbf N}$ are also bounded in the discrete setting, since
	\[  
	\left|A^{\mathbf\osi}_{\mathbf N}1_{N\times M}A^{\mathbf\osi}_{\mathbf N}(k,\mathbf l)\right| = \left|\sum_{{\mathbf i} \in I_{\mathbf N}} A(k,\mathbf i)\sum_{j \in I_M} A(j,\mathbf l)\right| \leq M^2\Pi_{i=1,\dots,d}(2\bso_i\osi_i+1)(\bso_i+1).
	\]
	Here, $M$ is the largest element of $A$, i.e., $M=\|\bspline\|_{L^\infty(\R)}$,  $\Pi_{i=1,\dots,d}(2\bso_i\osi_i+1)$ is the number of points in $P^{\mathbf\osi}_{\mathbf N}$ that are in the support of $\phi_{\mathbf N}$, and $\Pi_{i=1,\dots,d}(\bso_i+1)$ is the number of elements in $\Phi_{\mathbf{N}}$ that overlap with a given point. Again, the upper bound is independent on $N$.
	
	Next, substituting $Z^{\mathbf\bso,\epsilon}_{\mathbf N}$ for $Z_{\mathbf N}$ in equation~\eqref{eq:plungestructure2}, we obtain 
	\begin{align}\label{eq:plungestructure3}
	\left| \, \left( A  - A(Z^{\epsilon})^*A \right)(\mathbf k,\mathbf l)  \, \right| = \left\{
	\begin{array}{ll}
	\mathcal O(\epsilon)& \text{if }\mysupp\phi_{\mathbf N,\mathbf l}\subset[0,1]^d\backslash\Omega,\\
	\mathcal O(\epsilon) & \text{if }\mysupp\phi_{\mathbf N,\mathbf l}\subset\Omega,\\
	|a_{\mathbf k,\mathbf l}^\epsilon| &\text{otherwise.}
	\end{array}\right.
	\end{align}
	with
	\[
	a^\epsilon_{\mathbf k,\mathbf l} = \sum_{\mathbf i\in I_{\mathbf N}} A(\mathbf k,\mathbf i)(\delta_{\mathbf i, \mathbf l} -(\phi^{\epsilon}_{\mathbf N,\mathbf l},\tilde\phi^{\epsilon}_{\mathbf N,\mathbf i})_{L^2(\Omega)} ) = a_{\mathbf k,\mathbf l} +\mathcal O(\epsilon).
	\]
	Here, all terms of size $\mathcal O(\epsilon)$ are bounded explicitly by $\epsilon C$ and in the discrete case they are bounded by $\epsilon \, C^{\mathbf\osi}$.
	
	For given $\mathbf k$ and $\mathbf l$, the $\mathbf i$th element in the sum for $a^\epsilon_{\mathbf k,\mathbf l}$ is non-zero if $A(\mathbf k,\mathbf i)$ is non-zero and if $\delta_{\mathbf i,\mathbf l}-(\phi^{\epsilon}_{\mathbf N,\mathbf l},\tilde\phi^{\epsilon}_{\mathbf N,\mathbf i})_{L^2(\Omega)}$ is non-zero. The former requires $\mathbf i$ to be in
	\[
	I_{1}(\mathbf k) = \{\mathbf j\in I_{\mathbf N} \, | \, \mysupp \phi_{\mathbf N,\mathbf k} \cap\mysupp \phi_{\mathbf N,\mathbf j}\neq\emptyset\}.
	\]
	The latter corresponds to $(\phi^{\epsilon}_{\mathbf N,\mathbf l},\tilde\phi^{\epsilon}_{\mathbf N,\mathbf i})_{L^2(\Omega)}$ being non-zero, which may only be the case if $\mathbf i$ is in
	\begin{align}
	I_2(\mathbf l) =\{\mathbf j\in I_{\mathbf N} \, | \, \mysupp \tilde \phi^{\epsilon}_{\mathbf N,\mathbf j}\cap \mysupp \phi_{\mathbf N,\mathbf l}\neq \emptyset \}.
	\end{align}
	The intersection $I(\mathbf k,\mathbf l)=I_1(\mathbf k) \cap I_2(\mathbf l)$  is non-empty if there exists a $\mathbf j$ such that
	\[
	\mysupp \phi_{\mathbf N,\mathbf k} \cap\mysupp \phi_{\mathbf N,\mathbf j}\neq\emptyset \quad \mbox{and} \quad \mysupp \tilde\phi^{\epsilon}_{\mathbf N,\mathbf j}\cap \mysupp\phi_{\mathbf N,\mathbf l}\neq \emptyset.
	\]
	This means, loosely speaking, that $\mathbf k$ is close to $\mathbf l$.
	
	Finally, we conclude that for a fixed $\mathbf l$ and by the above, $a^\epsilon_{\mathbf k,\mathbf l}$ is non-zero if $\mysupp\phi_{\mathbf N,\mathbf k}$ overlaps with the domain
	\begin{equation}
	\bigcup_{j\in I_{\mathbf N}\text{ with }  \mysupp \tilde\phi^{\epsilon}_{\mathbf N,\mathbf j}\cap \mysupp\phi_{\mathbf N,\mathbf l}\neq \emptyset}\mysupp\phi_{\mathbf N,\mathbf j}.
	\end{equation}
	This domain is a hypercube of which the length in the $d$th dimension grows like $\mathcal O(1/N_d)$. Since the support of all basis elements depends in the same way on $N$, the number of elements that overlap with this domain is independent of $N$  --- it does however depend on $\epsilon$ and $\mathbf\bso$. Therefore, also the number of non-zero elements $a^\epsilon_{\mathbf k,\mathbf l}$ for fixed $\mathbf  l$ is independent of $N$. Consequently, the number of elements in column $\mathbf l$ of $
	A-AZ^*A
	$
	that are in absolute value larger than $\delta\triangleeq\epsilon C$ is also independent of $N$. Since there are $\# K_{\mathbf N}(\Omega)$ non-zero columns in $A-AZ^*A$, the number of rows with elements larger than $\delta$ in absolute value is $\mathcal O(\# K_{\mathbf N}(\Omega))$. Also, the number of elements in $A-AZ^*A$ that are larger than $\delta$ is $\mathcal O(\# K_{\mathbf N}(\Omega))$.
	
	This concludes the full proof of the last two statements for the continuous case.
\end{proof}

The previous result can be adjusted for the case where an exact discrete dual is known with compact support such as in~\eqref{eq:compactsampleddiscrete}. We state it here as a corollary.

\begin{corollary}\label{col:sparse}
For the AZ pair  $(A^{\mathbf\osi}_{\mathbf N},Z^{\mathbf\osi}_{\mathbf N})$ where \[ (Z^{\mathbf\osi}_{\mathbf N})(k,\mathbf l) = \tilde\phi^{\mathbf\osi}_{\mathbf{N},\mathbf{l}}(\mathbf t_k)  \]
and the dual $\tilde\phi^{\mathbf\osi}$ has compact support, 
the matrix $A-AZ^*A$ has at most 
\begin{enumerate}
	\item $\# K_{\mathbf N}^{\mathbf\osi}(\Omega)$ non-zero columns
	\item $\RANK(A-AZ^*A) \leq \# K_{\mathbf N}^{\mathbf\osi}(\Omega)$
	\item $\mathcal O(\# K_{\mathbf N}^{\mathbf\osi}(\Omega))$ non-zero rows
	\item $\mathcal O(\# K_{\mathbf N}^{\mathbf\osi}(\Omega))$ non-zero elements.
\end{enumerate} 
The constant in the big $\mathcal O$ is independent of $N$.
\end{corollary}
\begin{proof}
	The corollary results from Theorem~\ref{thm:lowrank1} and using $\delta=0$ because of the compact support.
\end{proof}

\begin{assumption}\label{ass:boundary}
	In the remainder of this text we assume that the boundary $\partial \Omega$ of $\Omega \subset \mathbb{R}^d$ has dimension $d-1$, i.e., its dimension is one less than that of $\Omega$. With this we exclude, for example, fractal domains. This assumption allows the statement that the number of elements in both $\#K_{\mathbf N}(\Omega)$ and $\#K_{\mathbf N}^{\mathbf\osi}(\Omega)$ scales as $N^{\frac{d-1}{d}}$. In particular, the rank in the univariate case is constant as $N$ grows.
	
	Moreover, in the discrete case we assume a linear sampling rate $M = \mathcal O(N)$, with proportionality constant greater than $1$ to ensure oversampling. This allows asymptotic complexity estimates involving only the parameter $N$.
\end{assumption}

\begin{figure}
	\centering
	{\resizebox{.3\textwidth}{!}{\begin{tikzpicture}
\begin{axis}[ymode={log}, only marks]
    \addplot[mark={*}, mark color={black}]
        table[row sep={\\}]
        {
            \\
            1.0  4.473410137058122  \\
            2.0  4.473410137058101  \\
            3.0  1.044327151837626  \\
            4.0  1.0443271518376087  \\
            5.0  0.009960541560859418  \\
            6.0  0.009960541560852505  \\
            7.0  7.683266330765254e-14  \\
            8.0  7.19758845173239e-14  \\
            9.0  6.926518408648039e-14  \\
            10.0  6.896958390620621e-14  \\
            11.0  6.457066043475064e-14  \\
            12.0  6.450604877993427e-14  \\
            13.0  4.997514810070869e-14  \\
            14.0  3.9063895848334605e-14  \\
            15.0  3.7765512712448607e-14  \\
            16.0  3.5849404769841387e-14  \\
            17.0  3.557768508327276e-14  \\
            18.0  3.2934140712493494e-14  \\
            19.0  3.2450787799584214e-14  \\
            20.0  3.218743252127842e-14  \\
            21.0  3.205566802738634e-14  \\
            22.0  3.09723130364094e-14  \\
            23.0  2.631122668609889e-14  \\
            24.0  2.6242205537201564e-14  \\
            25.0  2.6000698803222737e-14  \\
            26.0  2.502291843276906e-14  \\
            27.0  2.3519306042000374e-14  \\
            28.0  2.3084047671694574e-14  \\
            29.0  2.286487632820455e-14  \\
            30.0  1.9111538164612603e-14  \\
            31.0  1.894272017033212e-14  \\
            32.0  1.8323061877412863e-14  \\
            33.0  1.7786084741636764e-14  \\
            34.0  1.7536600343752078e-14  \\
            35.0  1.731348826136467e-14  \\
            36.0  1.5917877665718017e-14  \\
            37.0  1.574704949792775e-14  \\
            38.0  1.4898829460231573e-14  \\
            39.0  1.4554438379590207e-14  \\
            40.0  1.3066770215598322e-14  \\
            41.0  1.3057454395706515e-14  \\
            42.0  1.2621712169751873e-14  \\
            43.0  1.2072759795973805e-14  \\
            44.0  1.2059495013125365e-14  \\
            45.0  1.2049062764902333e-14  \\
            46.0  1.2032773253819694e-14  \\
            47.0  1.1717323070110554e-14  \\
            48.0  1.1512853832089993e-14  \\
            49.0  1.1489945167314059e-14  \\
            50.0  1.1382443888476938e-14  \\
            51.0  1.123693628151933e-14  \\
            52.0  1.0708647806643603e-14  \\
            53.0  1.062263735441282e-14  \\
            54.0  1.05748440593713e-14  \\
            55.0  1.054691060580638e-14  \\
            56.0  9.842702244394427e-15  \\
            57.0  9.815642941182871e-15  \\
            58.0  9.702021943137558e-15  \\
            59.0  9.675337167296173e-15  \\
            60.0  9.602940822741817e-15  \\
            61.0  9.533329703595298e-15  \\
            62.0  9.248827871584963e-15  \\
            63.0  9.234895914982109e-15  \\
            64.0  9.209027101030783e-15  \\
            65.0  9.143544833058158e-15  \\
            66.0  8.961636158841838e-15  \\
            67.0  8.946864109598481e-15  \\
            68.0  8.732284450587411e-15  \\
            69.0  8.680846467043807e-15  \\
            70.0  8.577779197901653e-15  \\
            71.0  8.342865162498016e-15  \\
            72.0  8.120456839253661e-15  \\
            73.0  7.971827404328995e-15  \\
            74.0  7.706250129410312e-15  \\
            75.0  7.675634313352576e-15  \\
            76.0  7.652637014328332e-15  \\
            77.0  7.608820824776298e-15  \\
            78.0  7.58124227217744e-15  \\
            79.0  7.404027770122905e-15  \\
            80.0  7.360703114969298e-15  \\
            81.0  7.342694255288274e-15  \\
            82.0  7.31296797125668e-15  \\
            83.0  7.217478088318006e-15  \\
            84.0  6.930055881640953e-15  \\
            85.0  6.834054263532914e-15  \\
            86.0  6.676998600888047e-15  \\
            87.0  6.652690964338374e-15  \\
            88.0  6.547330231170225e-15  \\
            89.0  6.533868673427875e-15  \\
            90.0  6.399070323770178e-15  \\
            91.0  6.340875527552433e-15  \\
            92.0  6.299738170835316e-15  \\
            93.0  6.184013027315265e-15  \\
            94.0  6.171682290830626e-15  \\
            95.0  6.046364437519872e-15  \\
            96.0  5.973233733048866e-15  \\
            97.0  5.908811656668686e-15  \\
            98.0  5.869246647508254e-15  \\
            99.0  5.577323147548223e-15  \\
            100.0  5.562690618827369e-15  \\
            101.0  5.524066094612953e-15  \\
            102.0  5.244148414086905e-15  \\
            103.0  5.192514569868855e-15  \\
            104.0  4.9813886087622154e-15  \\
            105.0  4.911083636054951e-15  \\
            106.0  4.849487657857189e-15  \\
            107.0  4.837694528847894e-15  \\
            108.0  4.384924432341281e-15  \\
            109.0  4.3467323781758285e-15  \\
            110.0  4.2082494421886755e-15  \\
            111.0  4.043472159070739e-15  \\
            112.0  3.962747567690718e-15  \\
            113.0  3.823719208102827e-15  \\
            114.0  3.7689781865700236e-15  \\
            115.0  3.664453350817507e-15  \\
            116.0  3.49617168532312e-15  \\
            117.0  2.771459709761413e-15  \\
            118.0  2.764948287938598e-15  \\
            119.0  2.560423479607181e-15  \\
            120.0  2.2932494168410598e-15  \\
            121.0  1.777716558263634e-15  \\
            122.0  1.590723764293925e-15  \\
            123.0  1.1608605032472915e-15  \\
            124.0  1.1059397106721787e-15  \\
            125.0  8.128734095233895e-16  \\
            126.0  4.982221774549814e-16  \\
            127.0  2.8312830777343293e-16  \\
            128.0  1.27331295663721e-16  \\
            129.0  7.434401412195843e-17  \\
            130.0  4.862779584615421e-17  \\
            131.0  3.0030370658058115e-17  \\
            132.0  2.5496623094598844e-17  \\
            133.0  2.030732685708924e-17  \\
            134.0  1.9555422513950224e-17  \\
            135.0  1.6112749691267805e-17  \\
            136.0  1.5743535984052684e-17  \\
            137.0  1.436445726721444e-17  \\
            138.0  1.2149337784229242e-17  \\
            139.0  9.946037039456402e-18  \\
            140.0  7.352762447814222e-18  \\
            141.0  6.3864239913432245e-18  \\
            142.0  5.3209279702202914e-18  \\
            143.0  4.751760383155915e-18  \\
            144.0  3.973685001924812e-18  \\
            145.0  3.0619400135827646e-18  \\
            146.0  2.3714907719097638e-18  \\
            147.0  2.2587980927062777e-18  \\
            148.0  1.712497526972645e-18  \\
            149.0  1.4397393138881314e-18  \\
            150.0  1.0523528290544552e-18  \\
            151.0  9.154705632593205e-19  \\
            152.0  7.311385804403651e-19  \\
            153.0  6.489406727622512e-19  \\
            154.0  5.659234838663076e-19  \\
            155.0  4.2914615515241184e-19  \\
            156.0  6.528798082000061e-30  \\
            157.0  5.446103686492248e-30  \\
            158.0  5.289147818422297e-30  \\
            159.0  4.9568868544509314e-30  \\
            160.0  4.442674645081909e-30  \\
            161.0  3.87273312867336e-30  \\
            162.0  3.6723780887238915e-30  \\
            163.0  3.404702819136766e-30  \\
            164.0  3.1368178733851804e-30  \\
            165.0  2.9313923106857826e-30  \\
            166.0  2.776658755061315e-30  \\
            167.0  2.6548236484085686e-30  \\
            168.0  2.438833192832165e-30  \\
            169.0  2.1722601433911024e-30  \\
            170.0  1.9202503622641467e-30  \\
            171.0  1.7908372384674193e-30  \\
            172.0  1.5407490404099876e-30  \\
            173.0  1.3382100666167615e-30  \\
            174.0  1.190120229278151e-30  \\
            175.0  1.0542976382424463e-30  \\
            176.0  7.394314657266933e-31  \\
            177.0  3.7043115594845057e-31  \\
            178.0  3.265544626831903e-31  \\
            179.0  1.2139985586645586e-31  \\
            180.0  0.0  \\
            181.0  0.0  \\
            182.0  0.0  \\
            183.0  0.0  \\
            184.0  0.0  \\
            185.0  0.0  \\
            186.0  0.0  \\
            187.0  0.0  \\
            188.0  0.0  \\
            189.0  0.0  \\
            190.0  0.0  \\
            191.0  0.0  \\
            192.0  0.0  \\
            193.0  0.0  \\
            194.0  0.0  \\
            195.0  0.0  \\
            196.0  0.0  \\
            197.0  0.0  \\
            198.0  0.0  \\
            199.0  0.0  \\
            200.0  0.0  \\
        }
        ;
\end{axis}
\end{tikzpicture}}}
	\hspace{1em}{\includegraphics[width=.15\columnwidth,trim={0em -1.6em 0em 0em},]{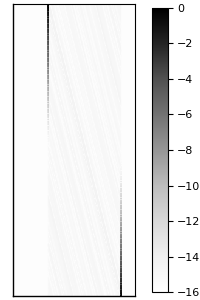}}
	{\resizebox{.3\textwidth}{!}{\input{lowrank2d.tikz}}}
	\hspace{1em}{\includegraphics[width=.15\columnwidth,trim={0em -1.6em 0em 0em},]{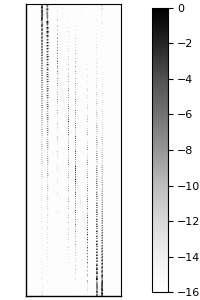}}
	\caption{\label{fig:sparseplunge}Panel 1,3: singular values of $A-AZ^*A$. Panel 2,4: logarithmic spy plot  of $|A-AZ^*A|(k,\mathbf l)$, $k\in I_M$, $\mathbf l=I_{\mathbf N}$. Panel 1,2: $\bso=3$, $N=200$, $\osi=2$, and $\Omega=[0.3,0.9]$. Panel 3,4: $\mathbf \bso=(3,3)$, $\mathbf N=(100,100)$, $\mathbf\osi=(2,2)$, and $\Omega$ is a disk with centre at $(0.5,0.5)$ and radius $1/3$.}
\end{figure}
The rank and sparsity structure of $A-AZ^*A$ are illustrated in Figure~\ref{fig:sparseplunge}.
The first panel shows the singular values of $A-AZ^*A$ for a 1-D domain ($N=200$). A small number of them are numerically significant, followed by a plateau of singular values near machine precision. By Theorem~\ref{thm:lowrank1}, in 1-D the number of numerically significant singular values remains constant with increasing $N$. The third panel shows a slightly larger number of large singular values for a 2-D domain ($N=100\times 100$), again followed by a drop to machine precision. Here, we expect $\mathcal O(N^{1/2})$ singular values to be above machine precision. The second and fourth panel illustrate the banded and sparse structure by displaying the spy plots of $A-AZ^*A$ in 1-D and 2-D. A pixel of light colour represent an matrix element with a value close to zero. 

In the following, we let 
$K$ equal $K^{\mathbf \bso}_N(\Omega)$ or $K^{\mathbf \bso,\mathbf \osi}_N(\Omega)$ in the continuous and discrete setting respectively.

\subsection{The sparse AZ algorithm}\label{ss:sparseaz}

The sparse AZ algorithm can be applied to the discrete problem if a compact dual is used to construct the matrix $Z$. Corollary~\ref{col:sparse} shows that the matrix contains only $\mathcal O(N^{(d-1)/d})$ non-zero elements, provided that Assumption~\ref{ass:boundary} is satisfied.   
Therefore, we consider an AZ algorithm where a sparse representation of $A-AZ^*A$ is created, which is solved using a direct least squares solver that relies on a rank-revealing sparse QR factorization~\cite{Davis2009}.

The algorithm can also be used by truncating the duals as in~\eqref{eq:truncateddual}. Then, the elements in $A-AZ^*A$ with absolute value smaller than some $\delta$ (recall $\delta\triangleeq \epsilon C$ in the proof of Theorem~\ref{thm:lowrank1}) in Theorem~\ref{thm:lowrank1} are zero. 

Following procedure allows to construct $A-AZ^*A$ as a sparse matrix with a computational complexity of $\mathcal O(N)$:
\begin{enumerate}
	\item Find $K$ in $\mathcal O(N)$ operations by iterating once over all basis functions and define $E$ as the $\#K\times N$ matrix that extends $K$ to $I_{\mathbf N}$.
	\item Find the index set
	\[
	\hat K= \{k\in I_M \, |\, l\in K\, :\,  AE(k,l) = 0 \}
	\]
	in~$\mathcal O\left(N\right)$ operations by iterating once over all $M$ points and define $R$ as the $\#\hat K\times M$ matrix that restricts $I_M$ to $\hat K$.
	\item Construct the sparse matrix $RAE$ in $\mathcal O\left(\# K\right)$.
	Using the structure of the supports of the bases involved, we know a priori the location of the non-zero elements in each of the $K$ columns and $\hat K$ rows.
	\item Find the indices
	\[ 
	\tilde K = \{ \mathbf l\in I_{\mathbf N}\, |\, k\in \hat K \, :\, Z(k,\mathbf l)=0   \}
	\] in~$\mathcal O\left(N§\right)$ operations and define $\tilde E$ as the $N\times\#\tilde K$ matrix that extends $\tilde K$ to $I_{\mathbf N}$. 
	\item Construct the sparse matrices $(RZ\tilde E)^*$ and $A\tilde E$ in $\mathcal O\left(\# K\right)$.
	\item Construct $E^*(I-Z^*A)E$ and use the result to construct $A(I-Z^*A)E$. Note that 
	\[\tilde E^*(I-Z^*A)E=\tilde E^*E-(RZ^*\tilde E)^*RAE\]
	and 
	\[ A(I-Z^*A)E= A\tilde EE^*(I-Z^*A)E,\]  
	see, e.g., the proof of Theorem~\ref{thm:lowrank1}. In general, a sparse matrix-matrix multiply has a lower bound of $N^2$ since the multiplication of two sparse matrices may result in a dense matrix. However, because of the structure in $RAE$, $A\tilde E$ and $(RZ\tilde E)^*$, $\tilde E(I-Z^*A)E$ and $A(I-Z^*A)E$ will have $\mathcal O\left(\# K\right)$ non-zero elements (see Corollary~\ref{col:sparse}). Moreover, the computation of each of these elements involves a number of operations that is independent of $N$. This step takes $\mathcal O\left(\# K\right)$.
	\item Construct the sparse matrix $A-AZ^*A$ as $[A(I-Z^*A)E]E^*$.
\end{enumerate}

We can make no precise statements about the complexity of the direct sparse QR solver, other than the observation that the cost in our application seems to scale roughly proportional to the number of non-zero entries.\footnote{We quote the reasoning in~\cite{Davis2009} here : \emph{A solver based on a sparse QR of $M=QR$ is subdivided in several steps. For \textit{ordering and symbolic analysis} no tight bounds are known for the time complexity. Experiments show that they scale roughly proportionally with the number of non-zeros in $M$ if \texttt{CHOLAMD} is used. However, exceptions exist. Next follows a \textit{symbolic factorization} that takes in practice $\mathcal O(|M|+|R|)$ number of steps where $|R|$ is the integers of steps needed to represent the multifrontal structure of $R$. This step can be a lot faster than the creation of $A^*A$. The \textit{numerical factorization} consists of a series of dense QR factorizations. The number and size depends on the earlier analysis.}} Furthermore, our solver relies on the efficient implementation of the product of sparse matrices as a sparse matrix. We have only shown that the result is indeed highly sparse, not how the multiplication algorithm should be implemented. We observe that the standard multiplication of sparse matrices in the current version of Julia (version 1.2, September 2019) yields optimal complexity in practice.

\subsection{The reduced AZ algorithm}\label{ss:reduced}
As it is difficult to predict the behaviour of the direct sparse QR solver used in the sparse AZ algorithm, we also introduce an algorithm that uses a traditional pivoted QR solver. This solver does not take advantage of the full sparsity of $A-AZ^*A$, as it simply uses dense matrices instead of sparse matrices.
We do exploit some sparsity of $A-AZ^*A$ by removing the zero columns and rows. In doing so, we obtain a system that has size $\mathcal O(N^{(d-1/d)}\times N^{(d-1)/d})$, see Corollary~\ref{col:sparse}.
To solve this system efficiently, it is not necessary to rely on randomized linear algebra since the rank is proportional to the dimension of the reduced matrix.

The algorithm is presented in Algorithm~\ref{alg:reducedAZ}.
\begin{algorithm}[t]
	\caption{The reduced AZ algorithm}\label{alg:reducedAZ}
	{\bf Input:} $A,Z \in \C^{M\times N}$, $b\in\C^M$\\
	{\bf Output:} $x\in\C^N$ such that $Ax \approx b$
	\begin{algorithmic}[1]
		\State Determine $R\in \{0,1\}^{\#I(K)\times M}, E\in \{0,1\}^{N\times \#K}$ such that $R(A-AZ^*A)E$ contains the non-zero columns and rows of $A-AZ^*A$
		\State Solve $R(I-AZ^*)AEx_1 = R(I-AZ^*)b$
		\State $x_2 \gets Z^*(b-AEx_1)$
		\State $x \gets Ex_1 + x_2$
	\end{algorithmic}
\end{algorithm}
In the first step, we determine the index set $K$. The selection of the non-zero columns is represented by the $N\times \#K$ matrix
\[
E(\mathbf k,\mathbf l)= \delta_{\mathbf k,\mathbf l}, \quad\mathbf k\in I_{\mathbf N},\mathbf l\in K.
\]
The index set $K$ can be constructed in ${\mathcal O}(N)$ operations most simply by iterating over all basis functions and checking their supports, and $E$ can be realized in practice as a sparse matrix.

After the non-zero columns, we also look for the indices of the non-zero rows,
\begin{equation*}
I(K)=\left\{k \in I_{M}\left| \, \forall \mathbf l\in K : \, \left|\left(I-AZ^*\right)AE\left(k,\mathbf l\right)\right|=0\right.\right\}.
\end{equation*}
This restriction operation is represented by the $\#I(K)\times M$ matrix, with $M$ the number of points in the discrete setting and $M=N$ in the continuous setting:
\[
R(k, l)= \delta_{k,l},\quad k\in I(K), \, l \in I_{M}.
\]
The computation of this index set $I(K)$ is more involved. It can be computed in $\mathcal O(N)$ operations using the sparse representation of $A-AZ^*A$ in the previous section or simply by checking the zero rows of the dense matrix $(A-AZ^*A)E$. The latter approach takes more operations, but this does not influence the complexity of the total number of operations of the reduced AZ algorithm as we also have to take into account the creation of the dense matrix and the solve step.

\begin{theorem}\label{thm:reducedaztimings}
	Provided Assumption \ref{ass:boundary} is satisfied, the reduced AZ algorithm (Algorithm~\ref{alg:reducedAZ}) using the AZ pairs $(A_{\mathbf N},Z_{\mathbf N})$ and $(A^{\mathbf\osi}_{\mathbf N},Z^{\mathbf\osi}_{\mathbf N})$ can be implemented with $\mathcal O(N^{3/2})$ operations in 2-D and $\mathcal(N^{3(d-1)/d})$ operations in $d$-D with $d>1$.
	In 1-D the number of operations depends on the dual used; for a dual computed by inverting the Gram matrix it is $\mathcal O(N{\log(N)})$; for a compact dual, it is $\mathcal O(N)$.  
\end{theorem}
\begin{proof}
	The extension matrix $E$ can be constructed in $\mathcal O(N)$ operations by iterating once over all basis functions. The matrix $R$ can be constructed most easily by explicitly forming the matrix $(A-A Z^*A)E$ of size $M \times n$, where $n=\# K = \mathcal O(N^{(d-1)/d})$ by Theorem~\ref{thm:lowrank1} and Assumption~\ref{ass:boundary}. The cost of constructing $(A- AZ^*A)E$ is $\mathcal O(n \texttt T_\text{mult})$. Selecting non-zero rows requires an additional $\mathcal O(nM)$ operations by iterating over all elements of the matrix. There are $m=\# I(K)= \mathcal O(N^{(d-1)/d}) = \mathcal O(\# K)$ such rows, hence $R \in \mathbb{R}^{m \times M}$ and forming the matrix $R(A - AZ^*A)E$ explicitly requires $\mathcal O(Mn)$ operations.
	
	Next, solving the $m\times n$ linear system with a direct solver requires $\mathcal O(mn^2)=\mathcal O(n^3)$ operations. Taken together, the full time complexity of the algorithm is
	\[
	\mathcal O(N + n\texttt T_\text{mult} +  nM + n^3)  = \mathcal O(n\texttt T_\text{mult} + n^3)
	\]
	where we take linear oversampling $M=\mathcal O(N)$ into account.
	
	It remains to determine $\texttt T_\text{mult}$. The matrix-vector product with $A$ takes $\mathcal O(N)$ operations since it is sparse. Indeed, in the continuous approach, each basis function $\phi_{\mathbf N,\mathbf k}$ has a support that overlaps with the support of a fixed number of other basis functions and, due to the scaling of splines with $N$, this number is independent of $N$. In the discrete approach, each basis function similarly is non-zero in a finite number of grid points independently of $N$. Thus, each column of $A_{\mathbf N} \in \mathbb{C}^{N \times N}$ and $A^{\mathbf\osi}_{\mathbf N} \in \mathbb{C}^{M \times N}$ has a fixed number of non-zero entries and there are $N$ columns. 
	
	In the case of a dual computed by inverting a Gram matrix (see~\S\ref{ss:efficient_fft}) the matrix $Z^*$ can, analogously be written  as the combination of a sparse and a circulant matrix.
	The multiplication by a circulant matrix is expedited using the FFT. The matrix-vector product with $Z^*$ can thus be performed in $\mathcal O(N\log(N))$ operations.
	
	Using the discrete compact dual, $Z^*$ can be created in $\mathcal O(1)$ (see~\S\ref{ss:efficient_fft}) and applied in $\mathcal O(N)$ operations instead of $\mathcal O(N\log(N))$.
\end{proof}

\begin{figure}
	\centering
	\resizebox{\textwidth}{!}{\input{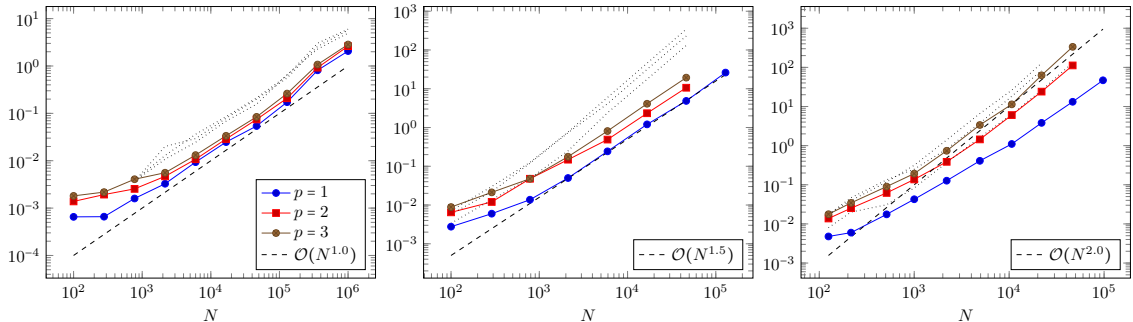}}
	\caption{
		\label{fig:AZR2timings1d-3d}
		Timings (in seconds) of the reduced AZ algorithm (Algorithm~\ref{alg:reducedAZ}) applied to the app\-rox\-imation problems of (left) $f(x)=e^x$ on $[0,1/2]$ using a 1-D spline basis on $[0,1]$, (middle) $f(x,y)=e^{xy}$ on $[0,1/2]^2$ using a 2-D spline basis on $[0,1]^2$, (right) $f(x,y,z)=e^{xyz}$ on $[0,1/2]^3$ using a 3-D spline basis on $[0,1]^3$. We approximate using splines of order $\bso=1$ (blue dots), $\bso=2$ (red squares) and $\bso=3$ (brown crossed dots). The expected asymptotic results of Theorem~\ref{thm:reducedaztimings} using the compact dual of to construct $Z$ are shown by the black dashed line: $\mathcal O(N)$ in 1-D, $\mathcal O(N^{3/2})$ in 2-D and $\mathcal O(N^{2})$ in 3-D. The  timings of using a non-compact (Gram-based) dual are present in the black dotted lines.}
\end{figure}
The time complexity of Theorem~\ref{thm:reducedaztimings} with AZ pair $(A^{\mathbf\osi}_{\mathbf N},Z^{\mathbf\osi}_{\mathbf N})$ is numerically illustrated in Figure~\ref{fig:AZR2timings1d-3d}. Shown are the timings of a 1-D, 2-D and 3-D spline extension approximation for increasing $N$. The figure compares the difference in using a compact and truncated Gram-based dual. While the time complexity results promised by the theorem are equal for both flavours of dual bases (except for the logarithmic factor in 1-D), we see that, in practice, the reduced AZ algorithm using the compactly supported basis takes less time. 

The computational cost is better than that of Fourier extension, which is $\mathcal O(N\log^2(N))$ in 1-D and $\mathcal O(N^2\log^2(N))$ in 2-D~\cite{matthysen2016fast,Matthysen2017}.

\section{Numerical experiments}\label{s:experiments}

All examples are performed using the \texttt{BSplineExtension.jl} package~\cite{coppe2019bsplineextension} written in Julia. Code producing the figures can also be found there. The continuous approach requires (multivariate) quadrature on irregular domains to evaluate the elements of $A_{\mathbf N}$ and $b_{\mathbf N}$. Therefore, we only consider experiments in the discrete setting here and in the following.

We compare the two algorithms proposed in \S\ref{s:algorithms} with each other and with the direct sparse QR solver.

\subsection{Time complexity}

\begin{figure}
	\centering
	\resizebox{\textwidth}{!}{\input{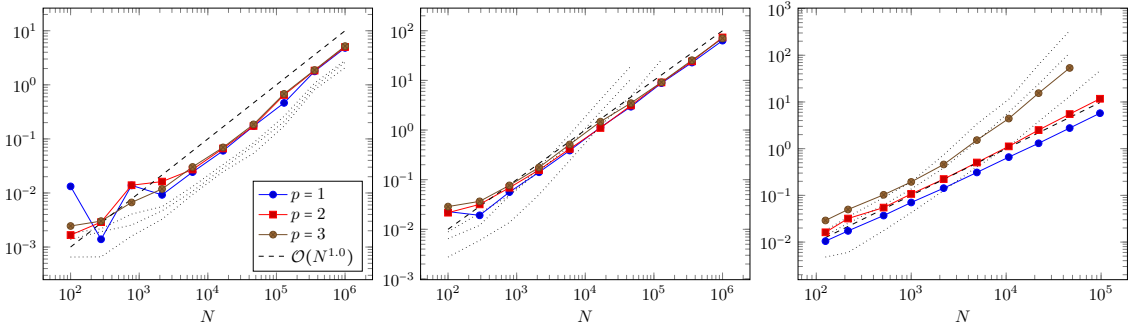}}
	\caption{\label{fig:AZStimings1d-3d}Timings (in seconds) of the sparse AZ algorithm  applied to the approximation problems of Figure~\ref{fig:AZR2timings1d-3d} and using the compact dual. In the black dashed line: $\mathcal O(N)$. The timings of Figure~\ref{fig:AZR2timings1d-3d} (using reduced AZ with a compact dual) are repeated in the black dotted lines. The sparse AZ algorithm is faster and seemingly exhibits linear complexity.}
\end{figure}

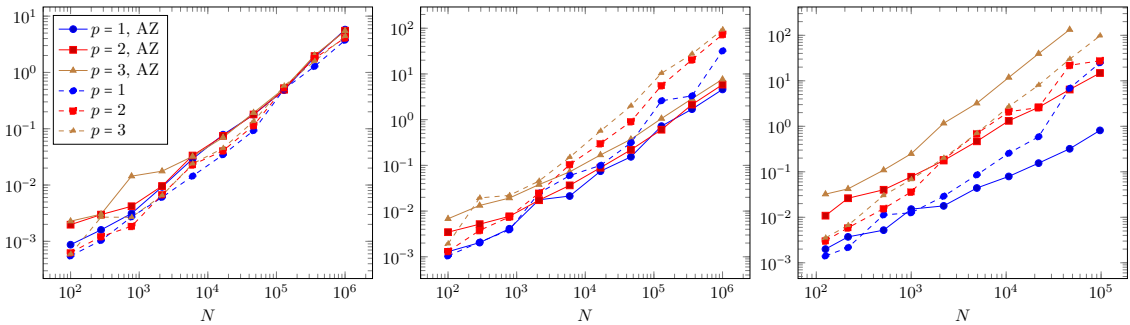
\begin{figure}
	\centering
	\resizebox{\textwidth}{!}{\begin{tikzpicture}
\begin{groupplot}[group style={group size={3 by 1}}]
    \nextgroupplot[xlabel={$N$}, xmode={log}, ymode={log}, legend cell align={left}, legend pos={north west}]
    \addplot
        table[row sep={\\}]
        {
            \\
            100.0  0.000870804  \\
            278.0  0.001595356  \\
            774.0  0.003122736  \\
            2154.0  0.009266708  \\
            5995.0  0.029081094  \\
            16681.0  0.078423852  \\
            46416.0  0.177974677  \\
            129155.0  0.483522939  \\
            359381.0  1.853666968  \\
            1.0e6  5.764175955  \\
        }
        ;
    \addplot
        table[row sep={\\}]
        {
            \\
            100.0  0.001975212  \\
            278.0  0.002981633  \\
            774.0  0.004196409  \\
            2154.0  0.009595867  \\
            5995.0  0.033423731  \\
            16681.0  0.074725971  \\
            46416.0  0.179197667  \\
            129155.0  0.493844178  \\
            359381.0  1.951830601  \\
            1.0e6  5.510899922  \\
        }
        ;
    \addplot[color=brown,mark=triangle*]
        table[row sep={\\}]
        {
            \\
            100.0  0.002281877  \\
            278.0  0.00303513  \\
            774.0  0.014376548  \\
            2154.0  0.017622304  \\
            5995.0  0.03222639  \\
            16681.0  0.069812269  \\
            46416.0  0.192220571  \\
            129155.0  0.519923546  \\
            359381.0  2.015400052  \\
            1.0e6  5.452452851  \\
        }
        ;
    \addplot[color=blue,mark=*,style=dashed]
        table[row sep={\\}]
        {
            \\
            100.0  0.000550625  \\
            278.0  0.001035136  \\
            774.0  0.002712595  \\
            2154.0  0.006044276  \\
            5995.0  0.014406626  \\
            16681.0  0.034535726  \\
            46416.0  0.092543114  \\
            129155.0  0.504742493  \\
            359381.0  1.279016366  \\
            1.0e6  3.730847775  \\
        }
        ;
    \addplot[color=red,mark=square*,style=dashed]
        table[row sep={\\}]
        {
            \\
            100.0  0.000628537  \\
            278.0  0.001212673  \\
            774.0  0.001850478  \\
            2154.0  0.006770587  \\
            5995.0  0.022888685  \\
            16681.0  0.040648897  \\
            46416.0  0.114036839  \\
            129155.0  0.535168731  \\
            359381.0  1.777809942  \\
            1.0e6  4.121840157  \\
        }
        ;
    \addplot[color=brown,mark=triangle*,style=dashed]
        table[row sep={\\}]
        {
            \\
            100.0  0.000593357  \\
            278.0  0.002685958  \\
            774.0  0.00267612  \\
            2154.0  0.00650511  \\
            5995.0  0.023838121  \\
            16681.0  0.045020564  \\
            46416.0  0.135531269  \\
            129155.0  0.576892168  \\
            359381.0  1.581233394  \\
            1.0e6  4.381494985  \\
        }
        ;
	\legend{{$p=1$, AZ},{$p=2$, AZ},{$p=3$, AZ},{$p=1$},{$p=2$},{$p=3$}}
    \nextgroupplot[xlabel={$N$}, xmode={log}, ymode={log}, legend cell align={left}, legend pos={north west}]
    \addplot
        table[row sep={\\}]
        {
            \\
            100.0  0.00129478  \\
            289.0  0.00204416  \\
            784.0  0.004105457  \\
            2116.0  0.017674844  \\
            5929.0  0.021334058  \\
            16641.0  0.074580897  \\
            46225.0  0.153754531  \\
            128881.0  0.722715843  \\
            358801.0  1.698174398  \\
            1.0e6  4.584304948  \\
        }
        ;
    \addplot
        table[row sep={\\}]
        {
            \\
            100.0  0.003458691  \\
            289.0  0.005177028  \\
            784.0  0.007686941  \\
            2116.0  0.017306559  \\
            5929.0  0.036568599  \\
            16641.0  0.089041834  \\
            46225.0  0.218926161  \\
            128881.0  0.60346034  \\
            358801.0  2.137525178  \\
            1.0e6  5.846569191  \\
        }
        ;
    \addplot[color=brown,mark=triangle*]
        table[row sep={\\}]
        {
            \\
            100.0  0.006790988  \\
            289.0  0.013272348  \\
            784.0  0.019293677  \\
            2116.0  0.037995724  \\
            5929.0  0.072198163  \\
            16641.0  0.168871416  \\
            46225.0  0.372134341  \\
            128881.0  1.058507727  \\
            358801.0  2.871499391  \\
            1.0e6  7.656833593  \\
        }
        ;
    \addplot[color=blue,mark=*,style=dashed]
        table[row sep={\\}]
        {
            \\
            100.0  0.001046181  \\
            289.0  0.002092634  \\
            784.0  0.00387968  \\
            2116.0  0.024761993  \\
            5929.0  0.059816552  \\
            16641.0  0.09995706  \\
            46225.0  0.311565707  \\
            128881.0  2.584265344  \\
            358801.0  3.317909796  \\
            1.0e6  31.872858894  \\
        }
        ;
    \addplot[color=red,mark=square*,style=dashed]
        table[row sep={\\}]
        {
            \\
            100.0  0.001316277  \\
            289.0  0.003779089  \\
            784.0  0.007451576  \\
            2116.0  0.024861499  \\
            5929.0  0.104345408  \\
            16641.0  0.299414094  \\
            46225.0  0.907515524  \\
            128881.0  5.563209395  \\
            358801.0  20.464771906  \\
            1.0e6  72.912459162  \\
        }
        ;
    \addplot[color=brown,mark=triangle*,style=dashed]
        table[row sep={\\}]
        {
            \\
            100.0  0.001926709  \\
            289.0  0.019429728  \\
            784.0  0.022048302  \\
            2116.0  0.04543762  \\
            5929.0  0.150571702  \\
            16641.0  0.563661934  \\
            46225.0  1.994102123  \\
            128881.0  10.363926781  \\
            358801.0  27.339271442  \\
            1.0e6  93.31786326  \\
        }
        ;
    \nextgroupplot[xlabel={$N$}, xmode={log}, ymode={log}, legend cell align={left}, legend pos={north west}]
    \addplot
        table[row sep={\\}]
        {
            \\
            125.0  0.001998877  \\
            216.0  0.003714894  \\
            512.0  0.005213411  \\
            1000.0  0.015115205  \\
            2197.0  0.017811831  \\
            4913.0  0.044055251  \\
            10648.0  0.079094207  \\
            21952.0  0.155046212  \\
            46656.0  0.320460571  \\
            97336.0  0.813853705  \\
        }
        ;
    \addplot
        table[row sep={\\}]
        {
            \\
            125.0  0.010859542  \\
            216.0  0.026241974  \\
            512.0  0.040251293  \\
            1000.0  0.077826767  \\
            2197.0  0.177016142  \\
            4913.0  0.468046671  \\
            10648.0  1.320093295  \\
            21952.0  2.618033423  \\
            46656.0  6.428339673  \\
            97336.0  14.851815537  \\
        }
        ;
    \addplot[color=brown,mark=triangle*]
        table[row sep={\\}]
        {
            \\
            125.0  0.032291317  \\
            216.0  0.042135786  \\
            512.0  0.109008052  \\
            1000.0  0.249813084  \\
            2197.0  1.166208274  \\
            4913.0  3.215704182  \\
            10648.0  11.841142289  \\
            21952.0  39.500359576  \\
            46656.0  133.509410193  \\
            97336.0  0.0  \\
        }
        ;
    \addplot[color=blue,mark=*,style=dashed]
        table[row sep={\\}]
        {
            \\
            125.0  0.001416217  \\
            216.0  0.002163298  \\
            512.0  0.011275416  \\
            1000.0  0.012563457  \\
            2197.0  0.029022277  \\
            4913.0  0.085653805  \\
            10648.0  0.255530912  \\
            21952.0  0.588860461  \\
            46656.0  6.864246876  \\
            97336.0  24.908394767  \\
        }
        ;
    \addplot[color=red,mark=square*,style=dashed]
        table[row sep={\\}]
        {
            \\
            125.0  0.00303972  \\
            216.0  0.005883818  \\
            512.0  0.015505743  \\
            1000.0  0.036134594  \\
            2197.0  0.185996272  \\
            4913.0  0.689466728  \\
            10648.0  2.091947489  \\
            21952.0  2.581987516  \\
            46656.0  21.916681431  \\
            97336.0  27.627688824  \\
        }
        ;
    \addplot[color=brown,mark=triangle*,style=dashed]
        table[row sep={\\}]
        {
            \\
            125.0  0.003523972  \\
            216.0  0.006651457  \\
            512.0  0.030165097  \\
            1000.0  0.067581562  \\
            2197.0  0.200244805  \\
            4913.0  0.70245067  \\
            10648.0  2.711254216  \\
            21952.0  8.030962377  \\
            46656.0  29.90104528  \\
            97336.0  97.576081477  \\
        }
        ;
\end{groupplot}
\end{tikzpicture}}
	\caption{Timings (in seconds) of the sparse AZ algorithm  using the compact dual (first three lines) and the direct sparse QR solver (last three lines) applied to: Left panel, 1-D (the interval $f(x)=e^x$ on $[0,1/2]$); Middle panel, 2-D ($f(x,y)=e^{xy}$ on a disk with centre $(1/2,1/2)$ and radius $0.4$); and, Right panel, 3-D ($f(x,y,z)=e^{xyz}$ on a ball with centre $(1/2,1/2,1/2)$ and radius $0.4$) using a B-spline basis on $[0,1]^d$.}\label{fig:AZSAStimings1d-3d}
\end{figure}

The time complexity of the sparse AZ algorithm is numerically illustrated in Figure~\ref{fig:AZStimings1d-3d} and Figure~\ref{fig:AZSAStimings1d-3d}. The former compares the sparse AZ algorithm with the reduced AZ algorithm. The latter compares with the direct sparse QR solver applied to the system $Ax=b$.

We see in Figure~\ref{fig:AZStimings1d-3d} that the timings for the sparse AZ algorithm and the reduced AZ algorithm are comparable in 1-D. However, for the higher-dimensional problems the sparse AZ algorithm is clearly faster. In particular, it seemingly reaches linear complexity. As mentioned above, it is not straightforward to prove this complexity.

We see in Figure~\ref{fig:AZSAStimings1d-3d} that for the 1-D approximation, there is no significant difference in timings between the sparse AZ algorithm and the direct sparse QR solver for $Ax=b$ (left panel). In 2-D, the AZ algorithm is clearly faster. In 3-D, the AZ approach is much faster than a direct solve for low spline orders. For higher spline orders, the splines have larger support and this directly impacts the computational cost of the algorithms. We see that the number of degrees of freedom needs to be sufficiently high before the sparse AZ algorithm outperforms the direct solve. For $\mathbf\bso=(3,3,3)$ and $\mathbf N=(46,46,46)$, e.g., the number of non-zero elements in $A$ and $A-AZ^*A$ is 13.363.264 and 28.112.684, respectively, while it is 419.601.512 and 239.763.184 for  $\mathbf N=(129,129,129)$.

\subsection{Approximation accuracy as a function of time}\label{ss:accuracy}

\begin{figure}
	\centering
	\resizebox{\textwidth}{!}{\input{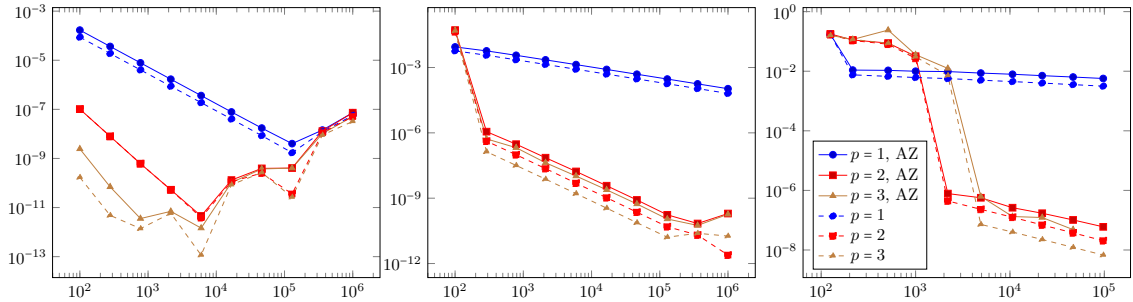}}
	\caption{Residuals of the systems solved in Figure~\ref{fig:AZSAStimings1d-3d}.}\label{fig:AZSASerrors1d-3d}
\end{figure}
For the experiment that compares the sparse AZ algorithm and the direct sparse QR approach, we also compare the residuals. Recall that the residuals relate directly to the approximation error. The residuals are comparable and show algebraic convergence in Figure~\ref{fig:AZSASerrors1d-3d}. The residuals of the AZ approach are slightly larger, consistently in all experiments.

\begin{figure}
	\centering
	{\includegraphics[width=.3\linewidth]{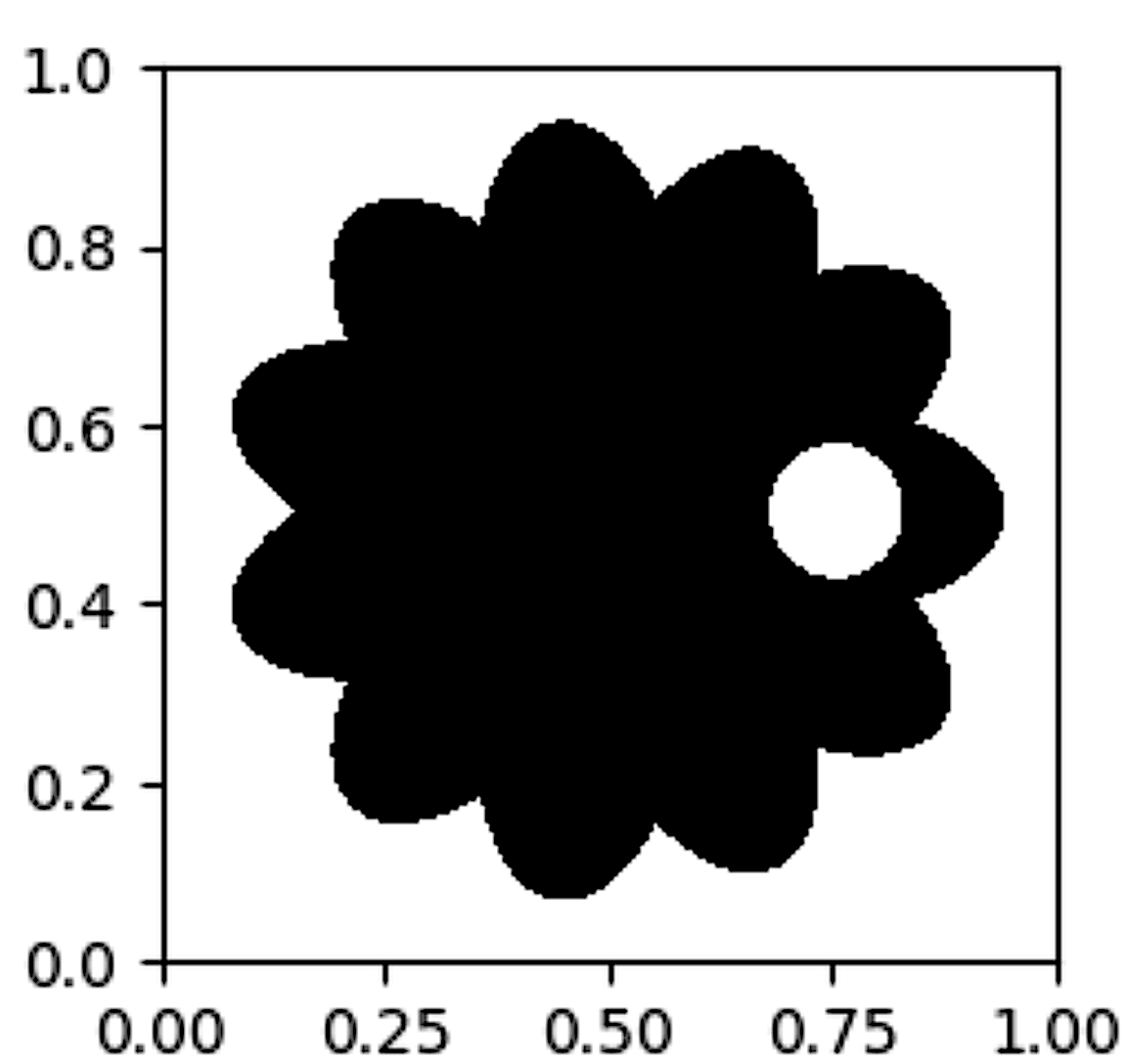}}
	\caption{The difference of the polar domain $\{x\in[-1,1]^2\, |\, \|x\|_2 \leq 0.35(2+0.5|\cos(5\arctan(x_2/x_1))|)\}$ with the disk with centre $(0.5,0)$ and radius $0.15$, affinely mapped from $[-1,1]^2$ to $[0,1]^2$. }\label{fig:flower}
\end{figure}

\begin{figure}
	\centering
	\resizebox{.6\textwidth}{!}{\input{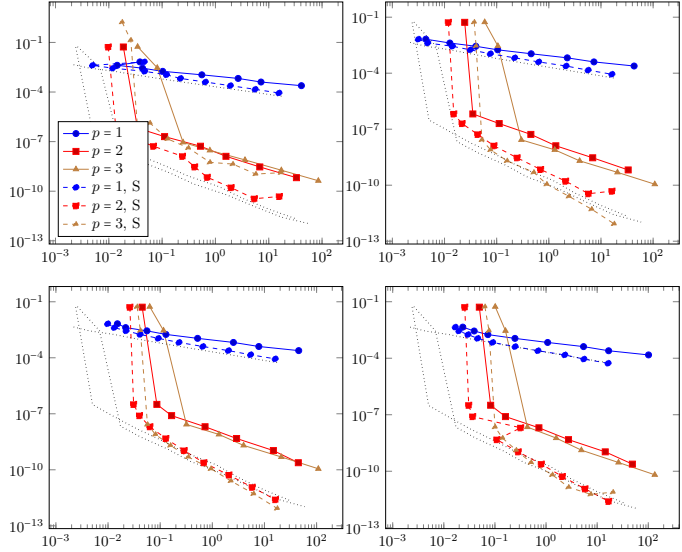}}
	\caption{Residuals as a function of timings of the reduced AZ algorithm using the compact dual and the sparse AZ algorithm (indicated by S) for spline extension approximation of $f(x,y)=e^{xy}$ on the flower-shaped domain of Figure~\ref{fig:flower}. In dotted lines, the direct sparse QR solver is used. Top left:  $M_{\text{abs}}=\infty$, top right: $M_{\text{abs}}=2$, bottom left: $M_{\text{abs}}=1$, bottom right:  $M_{\text{abs}}=0.8$.}\label{fig:AZcompactdualeff}
\end{figure}

Experiments indicate that the size of the residual is influenced by the choice of the compact dual spline basis, recall Remark~\ref{rem:dualsize}. This is illustrated in Figure~\ref{fig:AZcompactdualeff}. Theorem~\ref{thm:compactdual} guarantees that a dual spline exists with a certain minimal support. There may also be duals with larger support, but with smaller uniform norm of the discrete sequence, i.e., these are smaller duals. Here, we choose a dual with minimal support that also has uniform norm smaller than $M_{\text{rel}}=M_{\text{abs}}/\|b^\bso_\osi\|_{\infty}$, where $M_{\text{abs}}$ is a given constant. A lower uniform norm is beneficial for attaining smaller residuals. However, the larger support sizes lead to higher timings.

In order to illustrate the interplay between residual and clock time, we show in Figure~\ref{fig:AZcompactdualeff} the efficiency, i.e., the residual as a function of the timings. A lower curve indicates an improved efficiency. We study the effect of decreasing  $M_{\text{abs}}$ for the reduced and sparse AZ algorithm, and compare them with the efficiency of applying a direct sparse QR solver (which is independent of the choice of $M_{\text{abs}}$). 
The efficiency improves by requiring a smaller uniform norm. The change is most apparent in the sparse AZ algorithm, which for small $M_{\text{abs}}$ obtains the same efficiency as the direct sparse QR solver.

\begin{remark}
We have observed that the plain AZ algorithm (Algorithm~\ref{alg:AZ}) with a $Z$ matrix constructed using the compact dual basis performs very similarly to the reduced AZ algorithm in our implementation, even though the linear system in the former approach has (much) larger dimensions. This means that the algorithm based on randomized linear algebra to solve the low-rank problem may have implicitly removed the non-zero rows and columns automatically, resulting in an effective system of reduced size. This property of course simplifies the implementation of the scheme, but it depends on the implementation of the solver so it may not be a general property. We used the solver of the LowRankApprox.jl package in Julia~\cite{LowRankApprox.jl}.
\end{remark}

\section{Application to data smoothing on irregular geometries}\label{s:elevation}

\begin{figure}[t]
	\centering
	{\includegraphics[trim={1cm 0cm 1.5cm 1cm},clip,width=.4\linewidth]{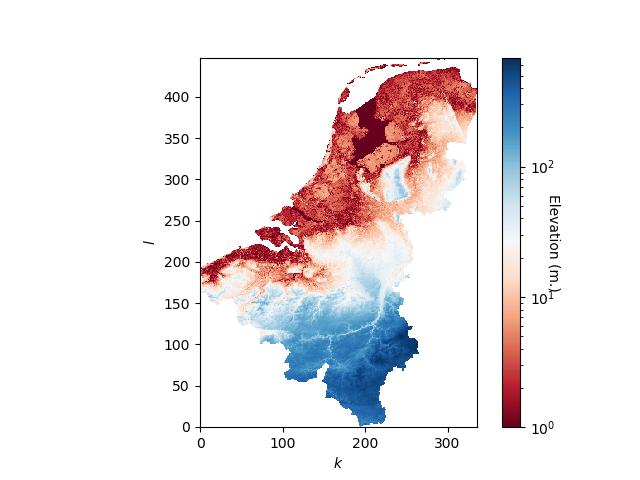}}
	\includegraphics[trim={1cm 0cm 1.5cm 1cm},clip,width=.4\linewidth]{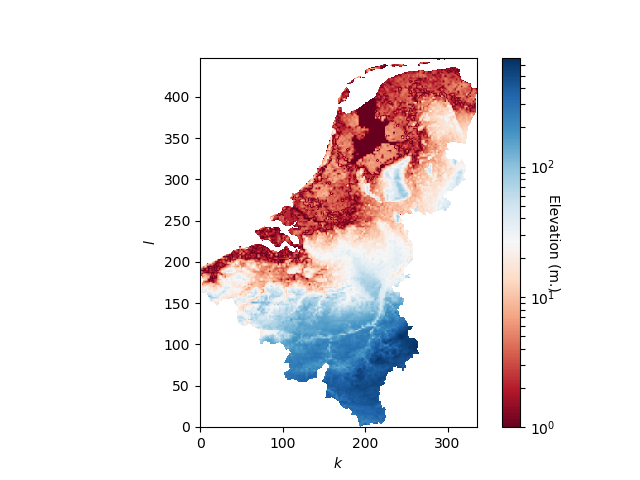}\\%
	\includegraphics[trim={1cm 0cm 1.5cm 1cm},clip,width=.4\linewidth]{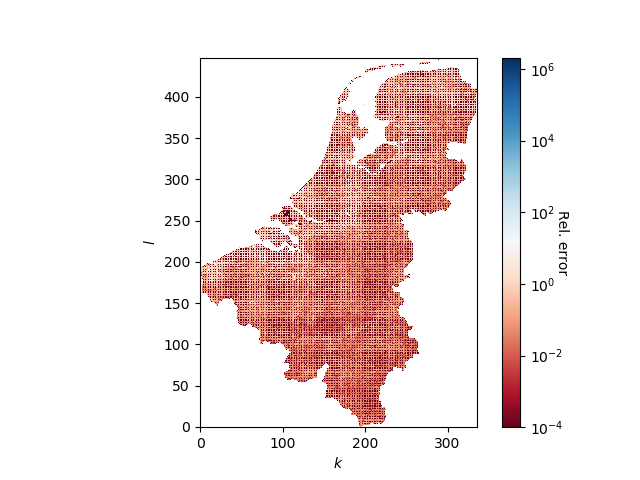}
	\includegraphics[trim={1cm 0cm 1.5cm 1cm},clip,width=.4\linewidth]{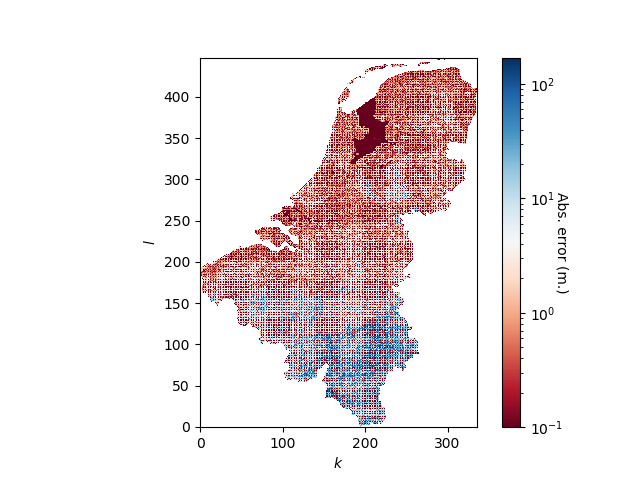}
	\caption{\label{fig:BENE}Top left: The elevation data of Belgium and The Netherlands. The data uses EPSG 3035/ ETRS89 as its coordinate system. The location data is equispaced. It contains the points that have coordinates $3799048.10 + 1000.73\, k$,  $k=1,\dots,336$, in the horizontal direction, and coordinates $2941371.64 + 999.85\, l$, $l=1,\dots,448$, in the vertical direction.
Top right: Approximation of the data using a 2-D spline extension with $\mathbf \bso=(1,1), \mathbf \osi=(2,2), \mathbf N=(168,224)$ and the sparse AZ algorithm. Elevation data smaller than 1 is truncated to 1 in the two figures at the top to make them more clear. Bottom left: The relative error. Bottom right: The absolute error. Relative errors smaller than $10^{-4}$ are truncated to $10^{-4}$ and absolute errors smaller than $10^{-1}$ are truncated to $10^{-1}$ to make the figure more clear.  The data was produced using Copernicus data and information funded by the European Union --- EU-DEM layers~\cite{BENEdata}.
	}
\end{figure}

We have invoked the final sparse AZ algorithm for the spline approximation of the elevation of the countries Belgium and The Netherlands. The data is provided by the Digital Elevation Model over Europe from the GMES RDA project (EU-DEM)~\cite{BENEdata}. The data consists of elevation (in meters) on an equispaced grid, with grid points confined by the highly irregular borders of these two countries. The results are depicted in Figure~\ref{fig:BENE}: shown  at the top left is a plot of the data, at the top right is the piecewise linear spline approximation. Shown at the bottom left and right is the relative and absolute error of the approximation in all the data points.

The approximation domain is not connected, due for example to a set of islands (de Waddeneilanden) belonging to The Netherlands in the north. Country borders are notoriously irregular, but so is the elevation data. Elevation data is more irregular in the hilly south of Belgium (the blue region covering the Ardennes) than in The Netherlands, which is mostly flat. This is adjusted for in the figure by showing the elevation in logarithmic scale, which reveals the low-elevation structure of the flat regions as well. Along the country borders on land, the regularity of the elevation data is not related to the regularity of the border, because the geography is of course continuous. In contrast, the two types of (ir)regularities are related in cases where the border is determined by a geographical feature, such as a sea (the North Sea) or a river. The dark homogeneous red region in the middle of the upper part of the figure has constant elevation because it is a lake at sea level (Ijsselmeer). It is included in the elevation data for The Netherlands because a large dike seals it off from the sea.

The least squares system in this example has dimensions $67970 \times 37632$, where $37632 = 168\times224$ is the number of degrees of freedom on a bounding box. The selection matrix $E$ is a $50176 \times 1386$ matrix. Thus, the system $(A - AZ^*A)E$ has dimensions $67970 \times 1368$, of which only $1736$ elements are non-zero. The time required to solve the final reduced linear system with the sparse AZ algorithm is only $0.22$ seconds on a MacBook Pro with a 2.8~GHz Quad-Core Intel Core i7 processor, 16~GB 1600 MHz memory and using Julia version 1.2. The relative residual of the solution is 0.0721. A direct sparse QR solver applied to the original linear system $Ax=b$ takes 0.48 seconds and results in a relative residual of 0.6537. The code for this experiment is available in~\cite{coppe2019bsplineextension}.

\section{Concluding remarks}

We have shown that spline extension frame approximation on general domains is possible and efficient using a regular spline basis defined on a bounding box. 
We used the AZ algorithm to isolate the ill-conditioned part in the computation of the approximation. This algorithm necessitates the knowledge of a dual basis to the spline basis. We compared different dual bases and their influence on the structure of $A-AZ^*A$ and conclude that an exact and compact dual is the best choice.  Spline extension leads to a low-rank and sparse structure of $A-AZ^*A$ that enables fast AZ-like algorithms.
The algorithm proposed in \S\ref{ss:reduced} has a provable complexity of $\mathcal O(N)$ operations in 1-D, $\mathcal O(N^{3/2})$ in 2-D and $\mathcal O(N^{3(d-1)/d})$ in $d$-D, $d>1$. This should be compared to the cubic cost of standard direct solvers. The other algorithms considered appear to offer better complexity, but this could not be shown.

The use of a direct sparse QR solver applied to the linear system $Ax=b$ without modification does not require a detailed study of the structure of the problem. As such, it is simpler to implement, and relevant to compare to. The experiment of Figure~\ref{fig:AZcompactdualeff} shows it to be more efficient, in the sense of achieving accuracy as a function of time. The small loss of accuracy due to the introduction of the compact dual basis is the main culprit of the decreased efficiency of the AZ-based algorithms. Indeed, the AZ-based method is faster time-wise. Other families of (more regular) dual discrete bases might find a better balance between accuracy and sparsity. 

We could use any basis on the bounding box. The aim of the follow-up work in \cite{waveletextension} is to use the techniques developed here to construct efficient algorithms using wavelets on the bounding box. This may pave the way for adaptive algorithms and refinement in the future.

Preliminary experiments seem to indicate that the solvers based on sparse QR factorizations loose their stability when using wavelets. The coefficient norm explodes and accuracy is lost. This may mean that the seemingly stable behaviour of the direct sparse QR solver in the experiments of this paper are specific to the regular structure of B-splines. Also, it is certainly not the case that a direct sparse QR solver can always solve a linear system at a cost that is linear in the number of non-zero entries. That this seems so for the case of spline extension is an observation that warrants further study.

\section*{Acknowledgements}
We benefited from useful discussions with Marcus Webb. We are also grateful to Eva Gadeyne who provided us with the elevation data for the example in~\S\ref{s:elevation}. The data was produced using Copernicus data and information funded by the European Union --- EU-DEM layers. 

\bibliographystyle{siam}
\bibliography{splineextension}

\begin{thebibliography}{10}

\bibitem{Adcock2019}
{\sc B.~Adcock and D.~Huybrechs}, {\em {Frames and Numerical Approximation}},
  SIAM Review, 61 (2019), pp.~443--473.

\bibitem{adcock2020frames}
\leavevmode\vrule height 2pt depth -1.6pt width 23pt, {\em Frames and numerical
  approximation {II}: generalized sampling}, arXiv e-prints,  (2020),
  p.~arXiv:1802.01950.

\bibitem{adcock2014stability}
{\sc B.~Adcock, D.~Huybrechs, and J.~{Mart{\'i}n-Vaquero}}, {\em On the
  numerical stability of {{Fourier}} extensions}, Found. Comput. Math., 14
  (2014), pp.~635--687.

\bibitem{Aldroubi1992}
{\sc A.~Aldroubi, M.~Unser, and M.~Eden}, {\em {Cardinal spline filters:
  Stability and convergence to the ideal sinc interpolator}}, Signal
  Processing, 28 (1992), pp.~127--138.

\bibitem{Averbuch2014}
{\sc A.~Z. Averbuch, P.~Neittaanm{\"{a}}ki, and V.~A. Zheludev}, {\em {Spline
  and Spline Wavelet Methods with Applications to Signal and Image Processing:
  Volume I: Periodic Splines}}, vol.~I, Springer, 2014.

\bibitem{Averbuch2016}
\leavevmode\vrule height 2pt depth -1.6pt width 23pt, {\em {Spline and Spline
  Wavelet Methods with Applications to Signal and Image Processing: Volume II:
  Non-Periodic splines}}, vol.~II, Springer, 2016.

\bibitem{Boyd2005}
{\sc J.~P. Boyd}, {\em {Fourier embedded domain methods: extending a function
  defined on an irregular region to a rectangle so that the extension is
  spatially periodic and C}}, Applied Mathematics and Computation, 161 (2005),
  pp.~591--597.

\bibitem{Bruno2007}
{\sc O.~P. Bruno, Y.~Han, and M.~M. Pohlman}, {\em {Accurate, high-order
  representation of complex three-dimensional surfaces via Fourier continuation
  analysis}}, Journal of Computational Physics, 227 (2007), pp.~1094--1125.

\bibitem{Chen1987a}
{\sc M.~Chen}, {\em {On the solution of circulant linear systems}}, SIAM
  Journal on Numerical Analysis, 24 (1987), pp.~668--683.

\bibitem{coppe2019bsplineextension}
{\sc V.~Coppe}, {\em {BSplineExtension.jl} v0.1.0}.
\newblock \url{https://github.com/FrameFunVC/BSplineExtension.jl}, 2019.

\bibitem{waveletextension}
{\sc V.~Copp\'e and D.~Huybrechs}, {\em Efficient function approximation on
  general bounded domains using wavelets on a cartesian grid}, arXiv
  preprint:2004.03537,  (2020, Submitted).

\bibitem{az}
{\sc V.~Copp\'e, D.~Huybrechs, R.~Matthysen, and M.~Webb}, {\em The {AZ}
  algorithm for least squares systems with a known incomplete generalized
  inverse}, SIAM Journal on Matrix Analysis and Applications,  (2020).
\newblock To appear.

\bibitem{Davis1979}
{\sc P.~J. Davis}, {\em {Circulant matrices}}, Pure and applied mathematics,
  Wiley, New York (N.Y.), 1979.

\bibitem{Davis2009}
{\sc T.~Davis}, {\em Multifrontral multithreaded rank-revealing sparse {QR}
  factorization}, in Dagstuhl Seminar Proceedings, Schloss
  Dagstuhl-Leibniz-Zentrum f{\"u}r Informatik, 2009.

\bibitem{de1976total}
{\sc C.~{De Boor}}, {\em {Total positivity of the spline collocation matrix}},
  Indiana University Mathematics Journal, 25 (1976), pp.~541--551.

\bibitem{demko1977inverses}
{\sc S.~Demko}, {\em Inverses of band matrices and local convergence of spline
  projections}, SIAM Journal on Numerical Analysis, 14 (1977), pp.~616--619.

\bibitem{devore1993constructive}
{\sc R.~A. DeVore and G.~G. Lorentz}, {\em Constructive approximation},
  vol.~303, Springer Science \& Business Media, 1993.

\bibitem{edelman1999}
{\sc A.~Edelman, P.~McCorquodale, and S.~Toledo}, {\em The {{Future Fast
  Fourier Transform}}?}, SIAM J. Sci. Comput., 20 (1999), pp.~1094--1114.

\bibitem{BENEdata}
{\sc {GMES RDA project (EU-DEM)}}, {\em Digital elevation model over europe
  (eu-dem)}, 2014.
\newblock data retrieved from European Environment Agency,
  \url{https://www.eea.europa.eu/data-and-maps/data/copernicus-land-monitoring-service-eu-dem}.

\bibitem{guide1978splines}
{\sc C.~d.~B. Guide}, {\em A practical guide to splines}, Applied Mathematical
  Sciences, 27 (1978).

\bibitem{hollig2003finite}
{\sc K.~H{\"o}llig}, {\em Finite element methods with B-splines}, SIAM, 2003.

\bibitem{hollig2001b}
{\sc K.~H{\"o}llig, U.~Reif, and J.~Wipper}, {\em {B}-spline approximation of
  {N}eumann problems}, Preprint, 2 (2001), p.~2001.

\bibitem{huybrechs2010fourierextension}
{\sc D.~Huybrechs}, {\em On the {F}ourier extension of non-periodic functions},
  SIAM J. Numer. Anal., 47 (2010), pp.~4326--4355.

\bibitem{LowRankApprox.jl}
{\sc JuliaMatrices}, {\em {LowRankApprox.jl} v0.4}.
\newblock \url{https://github.com/JuliaMatrices/LowRankApprox.jl}, 2019.

\bibitem{Lyon2011}
{\sc M.~Lyon}, {\em {A fast Algorithm for Fourier Continuation}}, SIAM Journal
  on Scientific Computing, 33 (2011), pp.~3241--3260.

\bibitem{matthysen2016fast}
{\sc R.~Matthysen and D.~Huybrechs}, {\em {Fast algorithms for the computation
  of Fourier extensions of arbitrary length}}, SIAM Journal on Scientific
  Computing, 38 (2016), pp.~A899----A922.

\bibitem{Matthysen2017}
\leavevmode\vrule height 2pt depth -1.6pt width 23pt, {\em Function
  approximation on arbitrary domains using {F}ourier extension frames}, SIAM
  Journal on Numerical Analysis, 56 (2018), pp.~1360--1385.

\bibitem{Parvizian2007}
{\sc J.~Parvizian, A.~D{\"{u}}ster, and E.~Rank}, {\em {Finite cell method}},
  Computational Mechanics, 41 (2007), pp.~121--133.

\bibitem{Schillinger2015}
{\sc D.~Schillinger and M.~Ruess}, {\em {The Finite Cell Method: A Review in
  the Context of Higher-Order Structural Analysis of CAD and Image-Based
  Geometric Models}}, Archives of Computational Methods in Engineering, 22
  (2015), pp.~391--455.

\bibitem{Schoenberg1969}
{\sc I.~J. Schoenberg}, {\em {Cardinal Interpolation and Spline Functions}},
  Journal of Approximation Theory, 2 (1969), pp.~167--206.

\bibitem{schumaker1981spline}
{\sc L.~Schumaker}, {\em Spline functions: basic theory. cambridge mathematical
  library}, 1981.

\bibitem{slepian1978}
{\sc D.~Slepian}, {\em Prolate spheroidal wave functions, {{Fourier}} analysis,
  and uncertainty {{V}}: The discrete case}, The Bell System Technical Journal,
  57 (1978), pp.~1371--1430.

\bibitem{Unser1991}
{\sc M.~Unser, A.~Aldroubi, and M.~Eden}, {\em {Fast B-Spline Transforms for
  Continuous Image Representation and Interpolation}}, IEEE Transactions on
  Pattern Analysis and Machine Intelligence, 13 (1991), pp.~277--285.

\bibitem{Unser1993}
\leavevmode\vrule height 2pt depth -1.6pt width 23pt, {\em {B-Spline Signal
  Processing: Part I - Theory}}, IEEE Transactions on Signal Processing, 41
  (1993), pp.~821--833.

\bibitem{Unser1993a}
\leavevmode\vrule height 2pt depth -1.6pt width 23pt, {\em {B-Spline Signal
  Processing: Part II - Efficient Design and Applications}}, IEEE Transactions
  on Signal Processing, 41 (1993), pp.~843--848.

\end{thebibliography}

\appendix
\section{Shift-invariant B-spline spaces}\label{app:A}
\subsection{Generators of spline spaces}

We refer the reader to~\cite{Averbuch2014,Averbuch2016} for an extensive description, in the setting of signal and image processing, of both periodic and non-periodic splines. 

The function $\phi(t)$ is called a generator for $S$ if every $f(t)\in S$ can be written as
\begin{equation*}\label{eq:continuousexpantion}
f(t)=\sum_{k\in\Z}c_k\phi\left(t-k\right).
\end{equation*}
There exist multiple generators for one space $S$. The B-spline $\bspline(t)$ is just one of them. Below we will introduce another which is biorthogonal to B-splines translated over integer shifts.

The $2\pi$-periodic function
\begin{equation*}\label{eq:charachteristic}
u(\omega) \triangleq \sum_{k\in\Z}e^{-i\omega k}\bspline(k)
\end{equation*}
is called the characteristic function of the space $S$ \cite[Definition (3.2)]{Averbuch2016}.
Because of B-spline symmetry and compact support, $u$ is a cosine polynomial. Examples are
\begin{equation*}
u^2(\omega) = 1,\quad 
u^4(\omega)=\frac{1+2\cos^2(\omega/2)}{3},\quad
u^5(\omega)=\frac{5+18\cos^2(\omega/2)+\cos^4(\omega/2)}{24}.
\end{equation*}
For our B-spline space, the characteristic function is strictly positive.

Another quantity that will be useful is the symbol. The symbol of a function $f$ is defined as
\begin{equation}\label{eq:defsymbol}
\xi(\omega)=\sum_{k\in\Z}e^{-i\omega k}c(k)
\end{equation}
where $c(k)$ are the coefficients of $f$ in the spline basis:
\begin{equation*}
f(t)=\sum_{k\in\Z}c(k)\bspline(t-k).
\end{equation*} The symbol of $\bspline(t)$ is 1.

\subsection{Exponential decay of dual B-splines}

We introduce the relation between the symbols of B-splines and their duals, and the characteristic function:
\begin{lemma}[{\cite[Proposition 4.7]{Averbuch2016}}]\label{thm:biorthogonal}
	Each generator $\phi(t)$ of $S$ has its dual counterpart $\tilde \phi(t)$ such that the biorthogonal relations
	\begin{equation}
	\int_{-\infty}^\infty \phi(t-k)\overline{\tilde\phi(t-l)}dt = \delta_{k,l}
	\end{equation}
	hold. The symbols $\tau(\omega)$ and $\tilde\tau(\omega)$ of the generators $\phi(t)$ and $\tilde\phi(t)$, respectively, are linked as
	\begin{equation}
	\tau(\omega)\overline{\tilde \tau(\omega)} = \frac{1}{u^{2\bso}(\omega)},
	\end{equation}
	where $u^{2\bso}(\omega)$ is the characteristic function of the space $S^{2\bso}$.
\end{lemma}

The following proof of Theorem~\ref{thm:expdecaycontdual} is not stated as a theorem in~\cite{Averbuch2016}, but it follows from the reasoning on page 60 which we repeat here in order to extend the reasonings further on.

\begin{proof}[Proof of Theorem~\ref{thm:expdecaycontdual}]
	By Lemma~\ref{thm:biorthogonal} the symbol $\tau(\omega)$ of $\tilde\bspline(t)$ is $\frac{1}{u^{2\bso}(\omega)}$. Note that its definition is~\eqref{eq:defsymbol}
	\begin{equation*}
	\tilde\tau(\omega)=\frac{1}{u^{2\bso}(\omega)} = \sum_{k\in\Z}c(k)e^{-i\omega k},
	\end{equation*}
	which shows that $c(k)$ are the Fourier coefficients of $\tilde\tau(\omega)$. The characteristic function $u^{2\bso}(\omega)$ is a $2\pi$ periodic cosine polynomial without zeros on $[0,2\pi]$ \cite[Lemma 6]{Schoenberg1969}. Therefore, $\tilde\tau(\omega)$ is analytic and periodic; and its Fourier coefficients decay exponentially as $|k|$ grows.
\end{proof}

\section{Dual B-splines with respect to discrete sampling}\label{app:B}

Since we are interested in discrete methods, we also study discrete sequences associated with sampled B-splines on a regular grid. Stability, fast algorithms and applications in this discrete setting were treated in~\cite{Unser1991,Unser1993,Unser1993a,Aldroubi1992}. These sampled B-splines are not to be confused with another set of discrete B-splines defined by discrete convolutions of rectangular pulses~\cite{Averbuch2016}. We consider duality with respect to a discrete inner product involving the sample points. For these dual B-splines we again show exponential decay in Theorem~\ref{thm:expdecaydiscretedual}.

\subsection{Notation}

In this section we adopt the notation of~\cite{Unser1991}.
The centered and shifted sampled B-spline sequences are defined by sampling the centered continuous B-spline with an integer oversampling factor $\osi$:
\begin{equation*}
b_\osi(k)=\bspline\left(k/\osi\right),\quad k\in\Z.
\end{equation*}
Sampling the continuous B-splines translated over integer shifts results in translates of the B-spline sequences by corresponding multiples of $\osi$,
\begin{equation}\label{eq:sampledshiftedsplines}
\bspline\left(\tfrac k\osi-l\right)=b_\osi(k-l\osi), \quad \forall k,l\in \Z.
\end{equation}

Like the B-splines, these discrete sequences have compact support. For $p>0$, we find from~\eqref{eq:support} that their discrete support is
\begin{equation}\label{eq:support_discrete}
\mysupp b_\osi = [-Q,Q],\qquad Q = \left\lceil \osi\frac{\bso+1}{2}-1\right\rceil.
\end{equation}
For $\bso=0$ the support is $\left\lceil -\tfrac{q}2 ,\tfrac{q-1}2  \right\rfloor$ which is not symmetric for even $\osi$.

In this text we use the convention that discrete convolution, defined by
\begin{equation*}
\left(a\star b\right)(k)=\sum_{l=-\infty}^\infty a(k-l)b(l),
\end{equation*}
takes precedence over evaluation {in order} to avoid a multitude of brackets and {to} ease notation. Thus $\left(a\star b\right)(k)=a\star b(k)$. If we introduce the shift operator $\delta_i(k)$ such that $\delta_i\star a(k)=a(k-i)$ we can write equation~\eqref{eq:sampledshiftedsplines} as
\begin{equation*}
\bspline\left(\tfrac k\osi-l\right) = \delta_{\osi l}\star b_\osi(k), \quad k,l\in\Z.
\end{equation*}

Next, we define upsampling by a factor of $\osi$ as
\begin{equation*}
[a]_{\uparrow \osi}(k) = \left\{\begin{array}{ll}
a(k')& k=\osi k'\\0&\text{otherwise,}
\end{array}\right.
\end{equation*}
and downsampling by the same factor as
\begin{equation*}
[a]_{\downarrow \osi}(k) = a(\osi k).
\end{equation*}

Analogously to the continuous case above, we can define a discrete shift-invariant space
\[
S_\osi=\myspan\{\Phi_\osi\}, \quad \mbox{with} \quad \Phi_\osi=\{\delta_{\osi k}\star b_\osi\}_{k\in\Z}.
\]
We call the sequence $g$ a generator for $S_\osi$ if every $f\in S_\osi$ can be written as
\begin{equation*}
f(k) = \sum_{l\in\Z}c(l)g(k-\osi l)
\end{equation*}
which in signal processing notation simplifies to
\begin{equation*}
f = [c]_{\uparrow \osi}\star g.
\end{equation*}
These spaces are invariant with respect to shifts by integer multiples of $\osi$.

\subsection{Discrete duality of sampled splines}

We are again interested in a dual generator; but this time based on a discrete inner product on the sequence space $S_\osi$,
\begin{equation}\label{eq:discrete_inner_product}
\llangle a, b \rrangle = \sum_{k\in\Z} a(k) \overline{b(k)}.
\end{equation}
We look for a dual generating sequence $\tilde b_\osi\in S_\osi$ that satisfies
\begin{equation}\label{eq:discretebiorthogonality}
\sum_{i=-\infty}^\infty \tilde b_\osi(i-\osi k) b_\osi(i-\osi l) =\llangle  \delta_{\osi k} \star \tilde b_\osi,  \delta_{\osi l}\star b_\osi\rrangle = \delta_{k,l}, \quad k,l\in\Z,{\quad \tilde b_\osi\in S_\osi}.
\end{equation}

Alternatively, we can describe the duality with respect to an inner product defined on the function space $S^\osi$. Since $\tilde b_\osi \in S_\osi$, the sequence corresponds to the samples of a continuous function $\tilde \bspline_\osi \in S^\osi$ that has a representation in the basis $\Phi$:
\begin{equation}\label{eq:discretespline_continuous}
\tilde \bspline_\osi(t) \triangleq \sum_{k\in\Z}c_\osi(k)\bspline(t-k).
\end{equation}
We obtain functions in $S^\osi$ that are biorthogonal in an oversampled equidistant grid:
\begin{equation*}
\langle\tilde \bspline_\osi(\cdot-k), \bspline_\osi(\cdot-l)\rangle_{\osi} = \delta_{k,l},
\end{equation*}
where we have used a discrete inner product\footnote{The discrete bilinear form defined here is clearly not an inner product on, say, $L^2(\mathbb{R})$, but it is on the shift invariant space $S$. That is because the discretization points include all integers, and the value of any function in $S$ is uniquely determined by its values at the integers due to the unique solvability of the cardinal B-spline interpolation problem. Hence, a spline that evaluates to zero at all integers is identically zero on the real line.}
\begin{equation}\label{eq:discreteinnerproduct}
\langle f, g\rangle_\osi=\sum_{k\in\Z}f\left(\tfrac k\osi\right)\overline{g\left(\tfrac k\osi\right)}.
\end{equation}

The $\osi$-shift biorthogonality property results in a reconstruction formula that is similar to the one expressed by~\eqref{eq:reconstruction} in the continuous case. We have that
\[
f(t) = \sum_{k \in \Z} \langle f, \tilde{\beta}_\osi\rangle_\osi \, \phi(t), \qquad f \in S^{\bso}.
\]
For this dual generator, the discretization of~\eqref{eq:discretespline_continuous} yields the expression
\begin{equation*}\label{eq:discreteexpansion}
\tilde b_\osi = [c_\osi]_{\uparrow \osi} \star b_\osi.
\end{equation*}
It turns out that the dual generator in $S_\osi$ is unique. We describe the solution in \S\ref{ss:discreteleastsquares} and \S\ref{ss:decay_discrete_dual}. We will consider other solutions of~\eqref{eq:discretebiorthogonality} later on in \S\ref{ss:compactdual}, by relaxing the requirement that the dual sequence lies in $S_\osi$.

\subsection{Discrete least squares approximations}\label{ss:discreteleastsquares}

In the interpolation problem for a function $f$ on the real line one wants to determine coefficients $y(k)$ such that
\begin{equation*}
f(k) = \sum_{l=-\infty}^\infty y(l)\bspline(k-l),\quad k\in\Z.
\end{equation*}
This problem is better known as the \emph{cardinal B-spline interpolation problem} and it is extensively investigated in~\cite{Schoenberg1969}. The interpolation problem is uniquely solvable.

We are interested in the solution of the more general, oversampled problem
\[
f(k/\osi) = \sum_{l=-\infty}^\infty y(l)\bspline(k/\osi-l),\quad k\in\Z
\]
for which {a solution} not necessarily exists {with equality to the samples of $f$}. Thus, we solve the problem in a least squares sense:
\begin{equation}\label{eq:lsproblem}
y = \argmin_{a}\sum_{k=-\infty}^\infty|f(k/q)-\sum_{l=-\infty}^\infty a(l)\bspline(k/\osi-l)|^2.
\end{equation}
This problem was investigated in~\cite{Unser1993}.

Written in a convolutional form we want to solve
\begin{equation}\label{eq:convlsproblem}
f_\osi(k) = b_\osi\star [y]_{\uparrow \osi}(k),
\end{equation}
where $f_\osi(k) = f(k/\osi)$, in a least squares sense.
In section IV.D of~\cite{Unser1993} {the solution $y(k)$ of the least squares problem~\eqref{eq:convlsproblem} is written in terms of filters}. We restate equations {\cite[(4.18)-(4.19)]{Unser1993}}:
\begin{equation}\label{eq:inversefiltersolution}
y(k) = s_\osi\star [b_\osi\star f_\osi]_{\downarrow \osi}(k),
\end{equation}
where
\begin{equation}\label{eq:snm}
s_\osi(k) = \left([b_\osi\star b_\osi]_{\downarrow \osi}\right)^{-1}(k),
\end{equation}
if the inverse of $[b_\osi\star b_\osi]_{\downarrow \osi}$ exists. The meaning of the inverse in this context is that the sequence $s_\osi$ satisfies
\begin{equation}\label{eq:inversesequence}
s_\osi \star [b_\osi\star b_\osi]_{\downarrow \osi}=\delta_0.
\end{equation}
The solution can be found by taking Fourier transforms
of both sides, writing the convolution as a product of terms, and then moving one factor to the other side by division. The Fourier transform of $s_\osi$ exists if the Fourier transform of $[b_\osi\star b_\osi]_{\downarrow \osi}$ does not vanish, else the division introduces poles. The latter Fourier transform is in fact strictly positive, as will be shown in the proof of Theorem~\ref{thm:expdecaydiscretedual} later on.

To that end, we define the Fourier transform of a sequence $f$ as $F(\omega) = \hat F(e^{i\omega})$, where $\hat F(z)$ is the $z$-transform of $f$:
\begin{align*}
\hat F(z) = \sum_{k\in\Z}f(k)z^{-k}.
\end{align*}

\subsection{Exponential decay of discrete dual B-splines}\label{ss:decay_discrete_dual}

We study the dual generators in $S_\osi$ which are $\osi$-shift biorthogonal to $b_\osi$ in more detail. First, we show that there is a generator that satisfies the conditions~\eqref{eq:discretebiorthogonality}. The coefficients $c(k)$ of this generator $\tilde b_\osi$ in $\Phi_\osi$ are obtained from the filter defined by~\eqref{eq:snm}.

\begin{theorem}\label{thm:discretedual}
	Provided that~\eqref{eq:inversesequence} is uniquely solvable for $s_\osi$, the generator $\tilde b_\osi(k)\in S_\osi$ that is $\osi$-shift biorthogonal to $b_\osi$ (with $\osi$-shift biorthogonality defined by~\eqref{eq:discretebiorthogonality}) is given by
	\begin{equation}\label{eq:discrete_dual_expansion}
	\tilde b_\osi(k)=[s_\osi]_{\uparrow \osi}\star b_\osi(k).
	\end{equation}
	Furthermore, the discrete B-spline least squares problem~\eqref{eq:lsproblem} is solved by applying $\tilde b_\osi(k)$ to $f_\osi$:
	\[
	y(k) = [\tilde b_\osi\star f_\osi]_{\downarrow \osi}(k).
	\]
\end{theorem}
\begin{proof}
	First we show that a sequence in $S_\osi$ of the form~\eqref{eq:discrete_dual_expansion} with coefficients equal to $s_{\osi}$ is $\osi$-shift biorthogonal to $b_\osi$.
	With the notation $a'(k) = a(-k)$ for time-reversal and $\overline{a}(k)=\overline{a(k)}$ for complex conjugation of a sequence, we can write the inner product~\eqref{eq:discrete_inner_product} as
	\begin{equation*}
	\llangle a, b\rrangle = \sum_{k=-\infty}^\infty a(k)\overline{b(k)} = a'\star \overline{b}(0)=a\star \overline{b}'(0).
	\end{equation*}
	Since $b_\osi(k)=b_\osi(-k)$ and $b_\osi(k)$ is real-valued, we can write for arbitrary $g=[d]_{\uparrow \osi}\star b_\osi\in S_\osi$
	\[\llangle  \delta_{\osi k} \star g,  \delta_{\osi l}\star b_\osi\rrangle = \delta_{\osi(k-l)} \star g \star \left(b_\osi\right)'(0)= g \star b_\osi(\osi(k-l))=\left[[d]_{\uparrow\osi}\star b_\osi\star b_\osi\right]_{\downarrow \osi}(k-l). \]
	If we use
	\begin{equation}\label{eq:downsampleofupsampled}
	[[a]_{\uparrow \osi}\star b]_{\downarrow \osi}(k) = \sum_{i\in\Z}b(\osi k-\osi i)a(i) = \sum_{l\in\Z}b(\osi l)a(k-l)=a\star[b]_{\downarrow \osi}(k)
	\end{equation}
	then we can simplify the expression: we find coefficients of a dual generator by solving
	\begin{equation*}
	d \star [b_\osi\star b_\osi]_{\downarrow \osi}=\delta_0
	\end{equation*}
	for $d(k)$. This is exactly the definition of $s_\osi$ given in~\eqref{eq:snm}.
	
	The second equality of the theorem is verified using the definition of $s_\osi$, ~\eqref{eq:downsampleofupsampled} and~\eqref{eq:inversefiltersolution}:
	\begin{equation*}
	[\tilde b_\osi\star f]_{\downarrow \osi}(k) = [[s_\osi]_{\uparrow \osi}\star b_\osi\star f]_{\downarrow \osi}(k) = s_\osi\star [b_\osi\star f]_{\downarrow \osi}(k)=y(k).
	\end{equation*}
\end{proof}

Knowing the form of the dual generator, we are ready to show that its coefficients $s_\osi$ decay exponentially.  To that end, we study its Fourier transform and those of $b_\osi$ and \begin{equation*}
c_\osi(k)=\bspline\left(k/\osi+\tfrac12\right).
\end{equation*} 
In the following lemma we first consider the Fourier transforms of $b_1$ and $c_1$. The lemma generalizes a result of \cite[Prop. 1]{Aldroubi1992} and is needed for the proof of the main theorem, Theorem~\ref{thm:expdecaydiscretedual}, further on.
\begin{lemma}\label{thm:sampledspline}
	Let
	\begin{align}
	B_1(\omega) &= \sum_{k\in\Z}b_1(k)e^{-i\omega k} = \sum_{k=-\infty}^\infty \left(\frac{\sin(\tfrac{\omega-2\pi k}2)}{\tfrac{\omega-2\pi k}2}\right)^{\bso+1},\\
	C_1(\omega) &=  \sum_{k\in\Z}c_1(k)e^{-i\omega k}= \sum_{k=-\infty}^\infty e^{i\tfrac{\omega-2\pi k}2}\left(\frac{\sin(\tfrac{\omega-2\pi k}2)}{\tfrac{\omega-2\pi k}2}\right)^{\bso+1}, \label{eq:Cn1}
	\end{align}
	be the Fourier transforms of $b_1(k)$ and $c_1(k)$ for $\bso\geq0$ and $\bso>0$, respectively.
	
	Then $B_1(\omega)$ and $D(\omega)=e^{-i\tfrac\omega2}C_1(\omega)$
	are strictly decreasing in $[0,\pi]$. Moreover, $B_1(\pi)>0$ for $p\geq 0$ and $D(\pi)=0$ for $p > 0$.
\end{lemma}
\begin{proof}
	Earlier results show that $B_1({\omega})$ is strictly positive and decreasing in $[0,\pi]$. See, e.g., \cite[{Lemma 6}]{Schoenberg1969}  and \cite[Prop. 1]{Aldroubi1992}. This only leaves us to show that $D(\omega)=e^{-i\tfrac\omega2}C_1(\omega)$ is strictly decreasing in $[0,\pi]$ and has a zero at $\omega=\pi$ for $\bso>0$.
	
	The continuous Fourier transform of $\bspline(t)$ is
	\begin{equation*}
	B(\omega) = \left(\frac{\sin(\tfrac{\omega}2)}{\tfrac{\omega}2}\right)^{\bso+1}
	\end{equation*}
	because of the recursive convolution. The sequences $b_1(k)$ and $c_1(k)$ are the sampled (and for $c_1(k)$, shifted over $\tfrac{1}{2}$) versions of $\bspline(t)$. Therefore, their Fourier transforms can be written as stated in the theorem.
	
	Noting that
	\[ \sin\left(\tfrac{\omega-2\pi k}{2}\right)=(-1)^k\sin\left(\tfrac{\omega}{2}\right) \]
	we simplify $D(\omega)$ to
	\[
	D(\omega)= \left(2\sin(\tfrac\omega 2)\right)^{\bso+1}\sum_{k=-\infty}^\infty (-1)^k(-1)^{k(\bso+1)}(\omega-2\pi k)^{-(\bso+1)}.
	\]
	If $\bso$ is even, $D(\pi)$ takes {the} following form:
	\begin{equation*}
	D(\pi) = \left(\frac{2}{\pi}\right)^{\bso+1}\sum_{k=-\infty}^\infty \frac{1}{(1-2 k)^{\bso+1}},
	\end{equation*}
	which is zero since the summation evaluates to zero:
	\begin{align*}
	\sum_{k=-\infty}^{-1} \frac{1}{(1-2 k)^{\bso+1}} = \sum_{k=1}^\infty \frac{1}{(1+2 k)^{\bso+1}} = \sum_{k=0}^\infty \frac{1}{(-1+2 k)^{\bso+1}}
	= \sum_{k=0}^\infty \frac{(-1)^{\bso+1}}{(1-2 k)^{\bso+1}}= -\sum_{k=0}^\infty \frac{1}{(1-2 k)^{\bso+1}}.
	\end{align*}
	For odd $\bso$ the reasoning is analogous. Therefore, $D(\omega)$ vanishes at $\omega=\pi$ for all $p>0$.
	
	Next, we expand $\frac{\partial D}{\partial \omega}(\omega)$:
	\begin{align*}
	\frac{\partial D}{\partial \omega}(\omega) &= (\bso+1)\sum_{k=-\infty}^\infty(-1)^k \left(\frac{\sin(\tfrac{\omega-2\pi k}2)}{\tfrac{\omega-2\pi k}2}\right)^{\bso}\frac{\tfrac12\cos(\tfrac{\omega}{2})(-1)^k\tfrac{\omega-2\pi k}{2}-\tfrac12\sin(\tfrac{\omega}{2})(-1)^k}{\left(\tfrac{\omega-2\pi k}2\right)^2}.
	\end{align*}
	The expression vanishes at $\omega=0$. For $\omega >0$, we split the sum into two parts:
	
	\begin{align}
	\frac2{\bso+1}\frac{\partial D}{\partial \omega}(\omega) &= (\cos(\tfrac{\omega}{2})\tfrac{\omega}{2}-\sin(\tfrac{\omega}{2}))\sum_{k=-\infty}^\infty \left(\frac{\sin(\tfrac{\omega-2\pi k}2)}{\tfrac{\omega-2\pi k}2}\right)^{\bso}\frac{1}{\left(\tfrac{\omega-2\pi k}2\right)^2}\nonumber \\
	&-\pi \sum_{k=-\infty}^\infty \left(\frac{\sin(\tfrac{\omega-2\pi k}2)}{\tfrac{\omega-2\pi k}2}\right)^{\bso}\frac{k}{\left(\tfrac{\omega-2\pi k}2\right)^2} \label{eq:twosums}
	\end{align}
	
	The factor in front of the first summation, $f(\omega)=\cos(\tfrac{\omega}{2})\tfrac{\omega}{2}-\sin(\tfrac{\omega}{2})$, is strictly negative for $\omega\in (0,\pi]$ since its derivative $f'(\omega)=-\tfrac14\sin(\tfrac\omega2)$ is strictly negative in $(0,\pi]$ and $f(0)=0$. Since both prefactors are negative, if both sums are shown to be positive we can conclude that $\frac{\partial D}{\partial \omega}(\omega)$ is negative, hence $D(\omega)$ is decreasing. That would conclude the proof.
	
	We consider the range $\omega \in (0,\pi]$. The first sum in~\eqref{eq:twosums},
	\begin{align*}
	S_1(\omega) = \sum_{k=-\infty}^\infty \left(\frac{\sin(\tfrac{\omega-2\pi k}2)}{\tfrac{\omega-2\pi k}2}\right)^{\bso}\frac{1}{\left(\tfrac{\omega-2\pi k}2\right)^2},
	\end{align*}
	has only positive terms when $\bso$ is even. For $\bso$ odd, we write $S_1(\omega) = \sum_{k=-\infty}^\infty (-1)^kf^p_k(\omega)$
	with
	\begin{align*}
	f^p_k(\omega)=	\left(\frac{\sin(\tfrac{\omega}2)}{\tfrac{\omega-2\pi k}2}\right)^{\bso}\frac{1}{\left(\tfrac{\omega-2\pi k}2\right)^2}.
	\end{align*}
	Because of the term $\frac12(\omega - 2\pi k)$ in the denominator, we find that $f^p_k(\omega)>0$ for $k \leq 0$ and $f^p_k(\omega)<0$ for $k>0$. Furthermore, due to the growth of the denominators of both factors for increasing $|k|$, it is true that
	\[ |f^p_{k-1}|<|f^p_{k}|, \quad\text{for }k\leq0,\qquad |f^p_{k+1}|<|f^p_{k}|, \quad\text{for }k>0.  \]
	We can then group terms in pairs with alternating signs. For $k > 0$:
	\begin{align*}
	(-1)^{2k-1}f^p_{2k-1}(\omega) + (-1)^{2k}f^p_{2k}(\omega) = -f^p_{2k-1}(\omega) +f^p_{2k}(\omega) = |f^p_{2k-1}(\omega)| -|f^p_{2k}(\omega)| > 0,
	\end{align*}
	and for $k \leq 0$:
	\begin{align*}
	(-1)^{2k-1}f^p_{2k-1}(\omega) + (-1)^{2k}f^p_{2k}(\omega) = -f^p_{2k-1}(\omega) +f^p_{2k}(\omega) = -|f^p_{2k-1}(\omega)| +|f^p_{2k}(\omega)| > 0.
	\end{align*}
	This shows that
	\begin{equation*}
	S_1(\omega) = \sum_{k=-\infty}^\infty(-1)^{2k-1}f^p_{2k-1}(\omega) + (-1)^{2k}f^p_{2k}(\omega) > 0
	\end{equation*}
	is also positive for odd $p$.
	
	For the second sum in~\eqref{eq:twosums}, we first note that $|\omega - 2k\pi| < |\omega  + 2k\pi|$ for $\omega\in(0,\pi]$, hence
	\begin{equation}\label{eq:mp}
	|f^p_k(\omega)| > |f^p_{-k}(\omega)|, \quad \text{for }k > 0.
	\end{equation}
	We also set out to show that
	\begin{equation}\label{eq:kf}
	(k+1)f^p_{k+1}(\omega)-kf^p_{k}(\omega)>0\quad\text{for } k>0,\qquad kf^p_k(\omega) - (k+1)f^p_{k+1}(\omega)>0\quad\text{for } k<-1.
	\end{equation}
	To that end, we start by rewriting
	\[
	kf^p_k(\omega) - (k+1)f^p_{k+1}(\omega) =
	\frac{k\left(1+\frac{1}{ k-\tfrac{\omega}{2\pi}}\right)^{p+2}-(k+1)}{(\omega-2\pi k-2\pi)^{p+2}}.
	\]
	Here, the denominator is negative for $k>0$ and positive for $k<-1$.
	For the numerator, we note that
	\[  \left(1+\frac{1}{k-f}\right)^{2}-\frac{k+1}{k} = -\frac{f^2-k(k+1)}{k(k-f)^2}  \]
	which is positive for $k>0$ and negative for $k<-1$, with $f=\omega/(2\pi)\in(0,1/2]$. Thus
	\[   \left(1+\frac{1}{k-f}\right)^{p+2}-\frac{k+1}{k} > \left(1+\frac{1}{k-f}\right)^{2}-\frac{k+1}{k}  > 0  \]
	for $k>0$ and
	\[   \left(1+\frac{1}{k-f}\right)^{p+2}-\frac{k+1}{k} < \left(1+\frac{1}{k-f}\right)^{2}-\frac{k+1}{k}  < 0  \]
	for $k<-1$. Multiplying with $k$ (and changing the direction of the inequality for negative $k$) shows that the numerator is positive and leads to the inequalities in~\eqref{eq:kf}.
	
	Returning to the second sum in~\eqref{eq:twosums}, it is positive for $p$ even since
	\[
	S^p_2(\omega) = \sum_{k=-\infty}^\infty kf^p_k(\omega) =\sum_{k=-\infty}^\infty k|f^p_k(\omega)| =\sum_{k=1}^\infty k\left(|f^p_k(\omega)|-|f^p_{-k}(\omega)|\right) >0,
	\]
	owing to~\eqref{eq:mp}. The sum is positive for $p$ odd since
	\begin{align*}
	S^p_2(\omega) &= \sum_{k=-\infty}^\infty kf^p_k(\omega) \\
	&=\sum_{k=1}^\infty \left(-\left(2k-1\right)f^p_{2k-1}(\omega)  + (2k)f^p_{2k-1}(\omega)\right) + \left((-2k)f^p_{-2k}(\omega)  - \left(-2k+1\right)f^p_{-2k+1}(\omega)\right)   \\
	&>0
	\end{align*}
	because of~\eqref{eq:kf}. This ends the proof.
\end{proof}

Previous work already proved the exponential decay of the dual generator coefficients of some cases. Earlier results include the dual taken with respect to the inner product $(f,g)$ defined in~\eqref{eq:L2innerproduct} see Theorem~\ref{thm:expdecaycontdual}; and, the dual with respect to integer grid, i.e., the discrete inner product~$\langle f,g\rangle_1$ of which the definition is given in~\eqref{eq:discreteinnerproduct}, see~\cite{Schoenberg1969}.
The exponential decay of the discrete dual generator coefficients in the oversampled case, $\langle f,g\rangle_{\osi}$ with~$\osi>1$, has, to the best of our knowledge, not been described in literature.

\begin{theorem}\label{thm:expdecaydiscretedual}
	For $\bso>0$,~\eqref{eq:inversesequence} is uniquely solvable for $s_\osi$ and the coefficients $s_\osi(k)$ decay exponentially as $|k|\to \infty$.
\end{theorem}
\begin{proof}[Proof of \ref{thm:expdecaydiscretedual}]
	The Fourier transform of~\eqref{eq:snm} follows from taking the Fourier transform of both sides in~\eqref{eq:inversesequence}. Noting that the latter involves the convolution of a sequence with itself, which results in squaring the corresponding Fourier transform. This is followed by downsampling by a factor of $\osi$, leading to \cite[eqn (4.20)]{Unser1993}, which we repeat here:
	\[
	S_\osi(\omega) = \left( \frac1\osi\sum_{k=0}^{\osi-1} B_\osi\left(\frac{\omega+2\pi k}{\osi}\right)^2\right)^{-1}.
	\]
	By Theorem~\ref{thm:discretedual} and since the Fourier coefficients of periodic functions analytic on the real line decay exponentially, the result holds if we can show that $S_\osi(\omega)$ is analytic on the interval $[0,2\pi]$ and by extension the real line, i.e., $1/S_\osi(\omega)$ has no zeros on the unit circle.
	
	To obtain an expression for $B_\osi$, we use known convolution expressions for the sampled B-splines~\cite{Unser1993}. They depend on the parity of $\bso$ and $\osi$. For $\osi$ odd:
	\begin{equation*}
	b_\osi(k)=\frac{1}{\osi}\underbrace{b^0_\osi\star b^0_\osi\star\dots \star b^0_\osi}_{\bso+1\text{ times}}\star b_1(k).
	\end{equation*}
	For $\bso$ odd and $\osi$ even:
	\begin{equation*}
	b_\osi(k)=\frac{1}{\osi}\delta_{(\bso+1)/2}\star\underbrace{b^0_\osi\star b^0_\osi\star\dots \star b^0_\osi}_{\bso+1\text{ times}}\star b_1(k).
	\end{equation*}
	Finally, for $\bso$ even and $\osi$ even:
	\begin{equation*}
	b_\osi(k)=\frac{1}{\osi}\delta_{(\bso+2)/2}\star\underbrace{b^0_\osi\star b^0_\osi\star\dots \star b^0_\osi}_{\bso+1\text{ times}}\star c_1(k).
	\end{equation*}
	
	The z-transform of the rectangular pulse is $\hat B^0_\osi(z)=z^{\lfloor \osi/2\rfloor}\left(\tfrac{1-z^{-\osi}}{1-z^{-1}}\right)$, where $\lfloor x\rfloor$ truncates to the smaller integer. Therefore, we arrive at the following expressions for $ B_\osi(\omega)=\hat B_\osi(e^{i\omega})$, which again depend on the parity of $\bso$ and $\osi$.
	When either $\bso$ or $\osi$ is odd:
	\begin{equation*}
	B_\osi(\omega) = \frac{1}{\osi} T_\osi(\omega) B_1(\omega).
	\end{equation*}
	Here, $T_\osi(\omega) = \hat T_\osi(e^{i\omega})$ with
	$\hat T_\osi(z) = z^{i_0}\left(\frac{1-z^{-\osi}}{1-z^{-1}}\right)^{\bso+1}$
	for $i_0=(\osi-1)(\bso+1)/2$ \cite[eqn (3.10)]{Unser1993}. This leads to
	\begin{equation*}
	T_\osi(\omega) = \left(\frac{\sin(\osi\tfrac\omega 2)}{\sin(\tfrac\omega 2)}\right)^{\bso+1}.
	\end{equation*}
	
	When both $\bso$ and $\osi$ are even:
	\begin{equation*}
	B_\osi(\omega) = \frac{e^{-i\omega/2}}{\osi}T_\osi(\omega)C_1(\omega).
	\end{equation*}
	
	The function $T_\osi(\omega)$ has zeros at $\omega = 2\pi \tfrac k\osi + 2\pi l$ for $k=1,\dots,\osi-1, l\in\Z$, while Lemma~\ref{thm:sampledspline} shows that $B_1(\omega)$ is strictly positive and $C_1(\omega)$ has zeros at $\omega = \pi+2\pi l$ with $l\in\Z$. Noting that the zero of $C_1$ coincides with one of the zeros of $T_\osi$ for $\bso$ even, we see that $B_\osi(\omega)$ has the same zeros as $T_\osi(\omega)$.
	
	Since \[K_{\osi,k}(\omega)\triangleq B_\osi\left(\frac{\omega+2\pi k}{\osi}\right)^2\] is a positive $2\pi \osi$-periodic function with zeros at $\omega=2\pi l$, $l\in\Z$, with the exception of $\omega=-2\pi k\osi +2\pi\osi l$, $l\in\Z$, $K_{\osi,k_1}$ and $K_{\osi,k_2}$ don't share zeros if $k_1\neq k_2$. Therefore, $L(\omega)\triangleq\tfrac1\osi\sum_{k=0}^{\osi-1}K_{\osi,k}(\omega)$ is $2\pi$-periodic and strictly positive. And, its inverse, $S_\osi(e^{i\omega})$, meets the previously stated requirements.
\end{proof}

\begin{theorem}\label{thm:completeness}
	Let $\bso=0$ and $\osi \neq 2$. Then~\eqref{eq:inversesequence} is uniquely solvable for $s_\osi$ and the coefficients $s_\osi(k)$ decay exponentially as $|k|\to \infty$. More specifically, for $\osi$ odd we have $s^0_\osi(k) = \tfrac1\osi\delta_0(k)$. For $\osi$ even, but $\osi \neq 2$, we have
	\[
	s^0_\osi(k) = \left\{
	\begin{array}{ll}
	\frac{-1}{(1-\osi)^{k+1}} & k \leq 0, \\
	0 & k > 0.
	\end{array}\right.
	\]
\end{theorem}
\begin{proof}[Proof of Theorem~\ref{thm:completeness}]
	If $\osi$ is even, then $b^0_\osi(k)=1$ for $-\osi/2\leq k\leq \osi/2-1$ and $0$ otherwise, see the comment under~\eqref{eq:support_discrete}. Hence, $b^0_\osi\star b^0_\osi(k)$ is non-zero for $-\osi \leq k\leq \osi-2$, and $[b^0_\osi\star b^0_\osi]_{\downarrow\osi}(k)$ is non-zero only if $-1 \leq k\leq 0$. To be precise,  $[b^0_\osi\star b^0_\osi]_{\downarrow\osi}(0)=\osi-1$ and $[b^0_\osi\star b^0_\osi]_{\downarrow\osi}(-1)=1$, therefore
	\begin{equation*}
	(S^0_\osi)^{-1}(\omega) = e^{-i\omega}+\osi-1.
	\end{equation*}
	The modulus is bounded above and below by $\osi-2\leq|(S^0_\osi)^{-1}(\omega)|\leq\osi$ for $\omega\in[0,2\pi]$. Thus, $(S^0_\osi)^{-1}(\omega)$ vanishes in $[0,2\pi]$ only if $\osi=2$ and $\omega=\pi$ at the same time. This means that $S^0_2(\omega)$ has a pole at $\omega=\pi$ and~\eqref{eq:inversesequence} is not uniquely solvable for $s^0_2$.
	
	Still for $\osi$ even, but $\osi \neq 2$, using the series expansion $\frac{1}{1-z}=\sum_{k=0}^\infty z^k$, $S^0_\osi$ can be written as
	\begin{align*}
	S^0_\osi(\omega) = \frac{1}{e^{-i\omega}+\osi-1} = -\sum_{k=0}^\infty\frac{e^{-i\omega k}}{(1-\osi)^{k+1}},
	\end{align*}
	so $|s^0_\osi(k)|$ decays exponentially for $k\rightarrow-\infty$ and $\osi\geq4$ and it is zero if $k>0$.
	
	If $\osi$ is odd, then $b^0_\osi(k)=1$ for $-(\osi-1)/2\leq k\leq (\osi-1)/2$. Hence, $b^0_\osi\star b^0_\osi(k)$ is non-zero for $-\osi+1 \leq k\leq \osi-1$, and $[b^0_\osi\star b^0_\osi]_{\downarrow\osi}(k)$ is non-zero if $k=0$. To be precise,  $[b^0_\osi\star b^0_\osi]_{\downarrow\osi}(0)=\osi$ and
	\begin{equation*}
	(S^0_\osi)^{-1}(\omega) = \osi.
	\end{equation*}
	Thus, $s^{0}_\osi(k)=\tfrac1\osi\delta_0(k)$.
\end{proof}

\subsection{Compact discrete duals to B-splines}\label{ss:compactdual}

We end this section with a description of duals to the sampled B-splines that have compact support. The $\osi$-shift biorthogonality conditions~\eqref{eq:discretebiorthogonality} admit alternative solutions outside of the shift-invariant space~$S_\osi$, in particular solutions with compact support. For these discrete solutions we are not able to provide a continuous analogue as in~\eqref{eq:discretespline_continuous}.

For clarity of the presentation we switch to a new notation for these duals. Thus, we are looking for a sequence $\tilde{h}_\osi$ that satisfies
\begin{equation}\label{eq:q-biorthogonality-v2}
\sum_{k=-\infty}^\infty \tilde h_\osi(k) b_\osi(k-\osi l) = \delta_l, \quad l\in\Z.
\end{equation}
If this holds, then $\osi$-shifts of $\tilde{h}_\osi$ define a dual generating sequence for a discrete shift-invariant space
\[
\tilde{U}_\osi = \SPAN \{ \delta_{\osi k} \star \tilde{h}_\osi\}_{k \in \Z}.
\]
This space may in general be different from $S_\osi$. We are encouraged by the observation that if both $\tilde h_\osi$ and $b_\osi$ are compact sequences,~\eqref{eq:q-biorthogonality-v2} reduces to a finite number of conditions.

\begin{proof}[Proof of Theorem~\ref{thm:compactdual}]
	Assume that $w$ is a sequence with support $[-K,K]$, i.e., that $w(k) = 0$ for $|k| > K$. Recall that the support of $b_\osi$ for $\bso>0$ is given by $[-Q,Q]$ with $Q= \lceil q \frac{\bso+1}{2}-1 \rceil$, see~\eqref{eq:support_discrete}. Substituting $w$ into~\eqref{eq:q-biorthogonality-v2} yields, with $L = \left \lfloor \frac{K+Q}{\osi} \right\rfloor$,
	\begin{equation}\label{eq:compact_system}
	\sum_{k=-K}^K w(k) b_\osi(k-\osi l) = \delta_0(l), \qquad -L \leq l \leq L,
	\end{equation}
	where the range of $l$ is restricted such that $|k-\osi l| > Q$ for $k \in [-K,K]$. Indeed, if $l > \frac{K+Q}{\osi}$, then
	\[
	k- \osi l < k - (K+Q) < -Q + (k-K) < -Q.
	\]
	The case $l < -\frac{K+Q}{\osi}$ similarly leads to $k - \osi l > Q$.  The conditions~\eqref{eq:compact_system} correspond to a linear system with $2L+1 = 2\left \lfloor \frac{K+Q}{\osi} \right\rfloor+1$ equations for $2K+1$ unknowns.
	
	We are interested in choosing $K$ such that there are more unknowns than conditions,
	\[
	2K + 1 \geq 2\left \lfloor \frac{K+Q}{\osi} \right\rfloor+1.
	\]
	This leads to $K \geq \left \lfloor \frac{K+Q}{\osi} \right\rfloor$ and hence $\frac{K+Q}{\osi} < K+1$, leading to
	\[
	K > \frac{Q-q}{q-1}.
	\]
	Substituting $Q= \lceil q \frac{\bso+1}{2}-1 \rceil$ results in
	\[
	K > \frac{q(p-1)-2}{2(q-1)} \quad\text{and}\quad K > \left\lceil\frac{q(p-1)-2}{2(q-1)}\right\rceil
	\]
	for odd and even $p$, respectively. This leads, after rearrangement, to the statements of the theorem.
	
	The system matrix of~\eqref{eq:compact_system} can also be written as
	\begin{equation*}
	A(l,k) = \beta\left(\tfrac{k}{\osi}-l\right),\qquad -L\leq l\leq L  ,\quad -K\leq k \leq K.
	\end{equation*}
	This is precisely the collocation matrix of a sequence of $2K+1$ spline functions $\bspline(\cdot -\tfrac{k}{q})$ evaluated in the $2L+1$ integers $-L\leq l\leq L$. The result follows if this system has full rank, i.e., if $A$ has rank $L$ (recall that $L < K$). We select a subset of $L$ splines, in such a way that each spline can be associated uniquely with one integer in its support. This is always possible because both the integers and the spline centers are regularly distributed on the same interval by construction. The corresponding $L \times L$ submatrix of $A$ represents a spline interpolation problem that is known to be uniquely solvable by \cite[Theorem 1]{de1976total}. Therefore, at least one exact solution of the underdetermined linear system exists.
\end{proof}

\begin{remark}
	Theorem~\ref{thm:compactdual} guarantees the existence of a dual sequence with small support, but does not state anything about its \emph{stability}, in the sense of having a bound on a discrete norm of $\tilde{h}_\osi$. In practice we observe that, owing to the norm-minimization property of least squares solutions, compact duals with smaller discrete norm can be found by solving~\eqref{eq:compact_system} in a least squares sense for a larger value of $K$ than the minimal one.
\end{remark}

\end{document}